\documentclass[11pt,leqno]{article}
\usepackage{authblk}
\usepackage{setspace} 
\usepackage[margin=1.2in]{geometry}
\usepackage{fancyhdr}
\usepackage{cite}
\usepackage{subfigure}

\usepackage{float}
\usepackage{graphicx}
\usepackage{graphbox}
\graphicspath{{figures/}} 
\usepackage{tikz} 
\usetikzlibrary{positioning, shapes.geometric}
\usetikzlibrary{arrows}

\usepackage{color}
\usepackage{xcolor}
\usepackage{hyperref}
\usepackage{indentfirst}
\usepackage{url}
\usepackage[]{cancel}

\DeclareGraphicsExtensions{.pdf,.png,.jpg}

\usepackage{amsmath,amssymb,amsthm,amsfonts}
\usepackage{algorithm}
\usepackage{algorithmic}
\newtheorem{theorem}{Theorem}[section]

\newtheorem{remark}[theorem]{Remark}
\newtheorem{example}[theorem]{Example}

\numberwithin{theorem}{section}
\numberwithin{equation}{section}
\numberwithin{figure}{section}
\numberwithin{algorithm}{section}
\usepackage{mathrsfs}
\usepackage{bm}
\usepackage{enumerate}
\makeatletter
\renewcommand{\p@subfigure}{}                    
\makeatother

\let\dotlessi\i  
\def\i{\mathrm{i}}
\def\ee{\mathrm{e}}

\def\d{{\mathrm d}}

\def\R {{\mathbb R}}

\def\le{\leqslant}
\def\ge{\geqslant}
\def\Omega{\varOmega}
\def\Delta{\varDelta}

\date{}

\begin{document}

\title{High-order mass-, energy- and momentum-conserving methods\\[2pt]
for the nonlinear Schr\"{o}dinger equation}

\author[1]{Georgios Akrivis}
\author[2]{Buyang Li}
\author[2]{Rong Tang}
\author[2]{Hui Zhang}

\affil[1]{\small Department of Computer Science and Engineering, University of Ioannina,
45110 Ioannina, Greece; and Institute of Applied and Computational Mathematics, FORTH,
70013 Heraklion, Crete, Greece.\ \ \texttt{akrivis@cse.uoi.gr}}
\affil[2]{\small Department of Applied Mathematics, The Hong Kong Polytechnic University,
Hung Hom, Hong Kong.\ \ \texttt{buyang.li@polyu.edu.hk}, \texttt{claire.tang@polyu.edu.hk},
\texttt{hui1203.zhang@polyu.edu.hk}}

\renewcommand\Authands{ and }
\date{}
\maketitle

\begin{abstract}

This paper introduces a novel formulation and an associated space-time finite element method for simulating solutions 
to the nonlinear Schr\"odinger equation. A major advantage of the proposed algorithm is its intrinsic ability to preserve 
the conservation of mass, energy, and momentum at the discrete level. This is proved for the numerical solutions determined 
by the fully discrete implicit scheme. An effective iterative scheme is proposed for solving the nonlinear system based on an equivalent 
formulation which suggests using Newton's iteration for the solution and no iteration for the Lagrange multipliers in the 
nonlinear system. Extensive numerical examples are provided to demonstrate the high-order convergence and effectiveness 
of the proposed algorithm in conserving mass, energy, and momentum in the simulation of one-dimensional Ma-solitons 
and bi-solitons, as well as of two-dimensional solitons governed by the nonlinear Schr\"odinger equation. The numerical results show that the mass-, energy- and momentum-conserving method designed in this paper also significantly reduces the errors of the numerical solutions in long-time simulations compared with methods which do not conserve these quantities. 

\end{abstract}


\bigskip


\section{Introduction}
This paper concerns the numerical solution of the nonlinear Schr\"{o}dinger (NLS) equation, subject to periodic boundary conditions, 
in a rectangular domain $\Omega \subset \mathbb{R}^d$ for some $d \ge1$. 
The equation can be expressed as 
\begin{subequations}
\begin{align}
\label{eeee1}
\i\partial_t u +\Delta u + f(|u|^2)u &= 0 \quad\hspace*{0.15cm} \text{ in } \Omega \times (0, T], \\
u|_{t=0} &= u_0 \quad \text{ in } \Omega,
\end{align}
\end{subequations}
where $u$ is an unknown complex-valued wave function, $\i = \sqrt{-1}$ denotes the imaginary unit, $u_0$ denotes the given 
initial value of the complex-valued wave function, and $f: \mathbb{R}^+ \rightarrow \mathbb{R}$ is a real-valued function, the derivative of a 
potential function $F: \mathbb{R}^+ \rightarrow \mathbb{R}$. Typical examples of the nonlinearity function $f$ take the form 
$f(|u|^2) = \mu|u|^{q-1}$, where $\mu \in \mathbb{R}$ and $q > 1$. The case where $\mu > 0$ is often referred to as the self-focusing 
model, in which the solution exhibits a finite-time blow-up when the initial energy is negative. Conversely, $\mu < 0$ is referred to as the self-defocusing model. 

The NLS equation \eqref{eeee1} is one of the fundamental equations in mathematical physics \cite{Bourgain,Ter2006,zabusky1965,pelinovsky1996,schurmann1996}. 
It is capable of describing the nonlinear dispersive waves in the modeling of the Bose--Einstein condensate \cite{bao2013,Schlein2010,Lieb2001}, the nonlinear optics 
\cite{chen2006,LeMesurier}, the deep-water modulation \cite{Peregrine,yuen1980}, and other applications. Accordingly, the numerical computation of the NLS 
equation has been extensively studied with various methods, including finite difference methods 
\cite{Akrivis93,Besse2004,Baoweizhu,Bao2004,Qinglin,Delfour81,Gao11,Wang13}, splitting methods \cite{Gauckler2010,Lubich2008,Thalhammer}, 
spectral methods \cite{Gong17}, discontinuous Galerkin methods \cite{Wenying2015,Xu2005,Guo2015}, and finite element methods 
(FEMs) \cite{KAD93,Caiwentao,Henning22,Karakashian,Karakashian99,Liuhai19,Tourigny91,Wang14,Zouraris01}.

It is widely recognized that solutions to the NLS equation possess conserved quantities such as mass, energy, momentum, etc., where the invariance of the momentum requires periodic boundary conditions.  In other words, 
for all $t \in (0, T]$, the following relations hold: 
\begin{equation}
M[u(t)] = M[u_0] ,\quad
E[u(t)] = E[u_0] \quad\text{and}\quad 
P[u(t)] = P[u_0] ,
\end{equation}
where 
\begin{align}
M[u] &= \frac{1}{2} \int_{\Omega} |u|^2 \, \d x ,\\
E[u] &= \frac{1}{2} \int_{\Omega} \left(|\nabla u|^2 - F(|u|^2)\right) \, \d x , \\
P[u(t)] &= \frac{1}{2} \text{Im} \int_{\Omega} \bar{u} \nabla u \, \d x ,
\end{align}
denote the mass, energy, and momentum of the solution, respectively. 
Correspondingly, the development of numerical methods that achieve high-order accuracy and retain the above-mentioned conservation properties at the discrete level simultaneously is important for long-time numerical simulation of wave propagation governed by the NLS equation. Let us mention that the periodicity is crucial for the invariance of the momentum. 

There has been extensive research on mass- and energy-conservative methods 
for the nonlinear Schr\"{o}dinger (NLS) equation. 
Typical mass-conserving methods include the implicit midpoint method \cite{hairer2000} and symplectic Rung--Kutta methods \cite{KAD93, Cooper1987, Hairer1993}. 
Additionally, energy conservation is commonly achieved using discrete gradient methods \cite{Gonzalez1996, McLachlan1999}, 
averaged vector field methods \cite{Quispel2008}, and continuous-stage Runge--Kutta methods \cite{Hairer2010}.

The modified Crank--Nicolson method \cite{Delfour81, Akrivis91, Henning17, sanz1984, Wang14} was one of the earliest and most widely 
adopted methods for preserving both mass and energy. To avoid solving nonlinear systems, a linearly implicit mass- and 
energy-conserving leap-frog scheme for NLS equations with cubic nonlinearity was developed in \cite{Fei1995}. This led to the 
design of linearly implicit mass- and energy-conserving methods, referred to as relaxation schemes \cite{Besse2004, Besse}, 
which exhibit second-order temporal convergence. The scalar auxiliary variable (SAV) technique was originally introduced 
to construct energy-decaying methods for dissipative equations \cite{Shen18, Shen19}, and was employed to develop 
schemes that conserve mass and SAV energy for Schr\"{o}dinger-type equations in \cite{Antoine2021, Feng21, Akrivis-DLi21, Guo2023}.

Bai et al.\ \cite{Bai24} developed high-order mass- and energy-conserving methods based on Gauss collocation in time 
and finite element discretization in space by correcting the numerical solution at every time level to conserve both mass 
and energy. Numerical results show that the solutions given by this prediction-correction procedure achieve higher accuracy 
in long-time simulations.

Ketcheson \cite{Ketcheson19} initially proposed a relaxation-type Runge--Kutta method to preserve quadratic energy, 
leading to the construction of numerous relaxation-type methods in 
\cite{Dongfang23, Dongfang223, Dongfang233, Ranocha20, Ranocha200, Ranocha2020}, typically preserving one invariant. 
Recently, new relaxation methods were developed in \cite{Biswas23} by combining the relaxation concept with embedded 
Runge--Kutta methods. Biswas and Ketcheson \cite{biswas2023} formulated an essentially explicit discretization for the NLS equation, conserving one or two invariants by integrating implicit-explicit higher-order Runge--Kutta time integrators 
with the relaxation technique and adaptive step size control. 

Despite these advancements, designing high-order numerical methods that conserve mass, energy, and momentum, or even more invariants of the NLS equation remains challenging and interesting. This paper aims to fill this gap by presenting a family of high-order methods for the NLS equation that conserve mass, energy, and momentum, based on novel formulations of the NLS equation leveraging Lagrange multipliers and associated constraint equations, as well as space-time FEMs that preserve the structure of the novel continuous formulation. The framework developed in this paper may also be further extended to conserve additional invariants of the NLS equation. The numerical results in this paper show that the mass-, energy- and momentum-conserving methods designed here also significantly reduce the errors of the numerical solutions in long-time simulations compared with methods which do not conserve these quantities; see Figure \ref{fig2-66}.

The paper is organized as follows: In Section \ref{sec2}, we introduce a novel formulation of the NLS equation 
utilizing Lagrange multipliers and a high-order space-time finite element algorithm for the novel formulation, 
as well as an iterative scheme for solving the  nonlinear system associated to the fully implicit numerical scheme. 
In Section \ref{sec4} we present several numerical examples that demonstrate 
the algorithm's high-order convergence and efficacy in conserving mass, energy, and momentum 
of the NLS equation, particularly in simulating Ma-solitons, bi-solitons, and solitons.

\section{The algorithm and main results}
\label{sec2}
In this section, we present the motivation and design of the algorithm that conserves mass, energy, and momentum of the NLS equation. 

\subsection{Motivation and continuous formulation}
Based on the mass conservation property, application of the fundamental theorem of calculus with respect to the time variable $t$ yields: 
\begin{equation}\label{mass-conservation}
  \begin{aligned}
  M[u(t_{n})] - M[u(t_{n-1})]  
  &= \int_{t_{n-1}}^{t_n} \frac{\d}{\d t}  M[u(t)]\, \d t \\
  &=  \int_{t_{n-1}}^{t_n}\text{Re} \int_{\Omega} \partial_t u \cdot \bar{u} \, \d x \d t \\
  &= - \int_{t_{n-1}}^{t_n} \text{Re} \int_{\Omega} \i\partial_t u \cdot \i\bar{u} \, \d x \d t 
  =0.
  \end{aligned}
\end{equation}
Similarly, the conservation properties of energy and momentum lead to 
\begin{equation}\label{energy-conservation}
 \begin{aligned}
   &E[u(t_{n})] - E[u(t_{n-1})]= \int_{t_{n-1}}^{t_n}  \frac{\d}{\d t}  E[u(t)]\, \d t\\
   &= -\text{Re} \int_{t_{n-1}}^{t_n} \int_{\Omega}\i\nabla\partial_t u \cdot \i\nabla\bar{u} \, \d x \d t 
   +\text{Re} \int_{t_{n-1}}^{t_n} \int_{\Omega}\i\partial_t u\cdot \i f(|u|^2)\bar{u}  \, \d x \d t\\
   &= \text{Re} \int_{t_{n-1}}^{t_n} \int_{\Omega}\i\partial_t u \cdot [\i\Delta\bar{u} + \i f(|u|^2)\bar{u}]  \, \d x \d t 
  =0
 \end{aligned}
\end{equation}
and, for the components $P_j[u]$, $j=1,\dots,d$, of the momentum,
\begin{equation}\label{momentum-conservation}
  \begin{aligned}
  P_j[u(t_{n})] - P_j[u(t_{n-1})] 
  &= \int_{t_{n-1}}^{t_n}   \frac{\d}{\d t}  P_j[u(t)] \,\d t \\
  &= - \int_{t_{n-1}}^{t_n} \text{Im} \int_{\Omega} \partial_t u \cdot \partial_j \bar u\, \d x \d t \\
  &= \int_{t_{n-1}}^{t_n} \text{Re} \int_{\Omega} \i\partial_t u \cdot \partial_j \bar u\, \d x \d t = 0.
  \end{aligned}
\end{equation}
Conversely, if the following constraints are satisfied:
\begin{subequations}\label{mass-momtentum}
\begin{align}
\text{Re}\int_{t_{n-1}}^{t_n} \int_{\Omega} \i\partial_t u \cdot \i\bar{u} \, \d x\d t &= 0 ,\\
\text{Re} \int_{t_{n-1}}^{t_n} \int_{\Omega}\i\partial_t u \cdot [\i\Delta\bar{u} + \i f(|u|^2)\bar{u}]  \, \d x \d t &=0 ,\\
\text{Re}\int_{t_{n-1}}^{t_n} \int_{\Omega} \i\partial_t u \cdot \partial_j \bar u\, \d x\d t &=0,
\end{align}
\end{subequations}
then the mass, energy, and momentum are conserved at the discrete time levels $t_m$. 


Our observation is that the constraints in \eqref{mass-momtentum} simply mean that $\i\partial_t u$ is orthogonal to the finite-dimensional subspace 
\begin{equation}
X(u)={\rm span}\{\i u, \i\Delta u+ \i f(|u|^2)u , \partial_1 u,\dotsc, \partial_d u\}  
\end{equation}
with respect to the space-time inner product of $L^2(\Omega\times(t_{n-1},t_n])$. 
Therefore, if we restrict the test functions for the NLS equation only to the orthogonal complement subspace $X(u)^\perp$, i.e., 
considering the following formulation of the NLS equation: 
\[
{\rm Re}\int_{t_{n-1}}^{t_n}\int_\Omega [\i\partial_t u + \Delta u + f(|u|^2)u] \bar v\, \d x\d t  = 0 \quad\text{for}\,\,\, v\in X(u)^\perp , 
\]
then this simply means that $\i\partial_t u + \Delta u + f(|u|^2)u \in X(u)$. Thus, there exist real-valued coefficients $\kappa_0,\kappa_1$ 
and $\kappa_{j+1}$, $j=1,\dotsc,d$, such that 
\begin{equation}\label{PDE-kappa}
\begin{split}
\i\partial_t u + \Delta u + f(|u|^2)u &= \kappa_0 \i u + \kappa_1 [\i\Delta u+ \i f(|u|^2)u]\\ 
&+ \sum_{j=1}^d \kappa_{j+1} \partial_j u \quad\text{in}\,\,\,\Omega\times(t_{n-1},t_n] . 
\end{split}
\end{equation}
The real-valued coefficients $\kappa_0,\kappa_1$ and $\kappa_{j+1}$, $j=1,\dotsc,d$, can be regarded as Lagrange multipliers for the constraints of the mass, energy, and momentum invariants. 

In the next subsection, we propose a fully discrete space-time FEM for the NLS equation based on the continuous formulations in \eqref{mass-momtentum} 
and \eqref{PDE-kappa}.

Throughout this article, we denote by $C$ a generic positive constant, possibly depending on the exact solution and $T$, but are independent of the mesh size and the time step size. The notation $X \lesssim Y$ means $X \leqslant C Y$ for some constant $C$.

\subsection{Space-time FEM}
%
For finite element discretization in space, we introduce a shape-regular and quasiuniform triangulation $\mathcal{T}_h$ of $\Omega$ with mesh size $h \in (0, 1]$. 
For any positive integer $r \geqslant 1$ and domain $K\subset\mathbb{R}^d$, $\mathbb{Q}^r(K)$ is  the space of complex-valued polynomials 
of degree up to $r$ in the domain $K$. Additionally, we denote by $S_h$ the periodic complex-valued Lagrange finite element space, i.e., 
\[
S_h = \{v \in C(\overline\Omega) : v|_K \in \mathbb{Q}^r(K)  \text{ for all } K \in \mathcal{T}_h \},
\]
where $C(\overline\Omega)$ denotes the space of uniformly continuous complex-valued functions in $\Omega$. 

We denote by $(\cdot,\cdot)$ and $\|\cdot\|$ the sesquilinear inner product and norm of the complex-valued Hilbert space $L^2(\Omega)$, i.e.,
\[ (u,v):=\int_\Omega u \overline{v}\, \d x \quad \textrm{and} \quad \|u\|=\sqrt{\int_\Omega |u|^2\, \d x }.\]
The discrete Laplacian operator on the finite element space $S_h$ is defined as the unique linear operator $\Delta_h: S_h \rightarrow S_h$ 
satisfying the following relation:
\begin{equation}
    \label{ee23}
 (\Delta_hv_h, w_h) = -(\nabla v_h, \nabla w_h) \quad \forall w_h \in S_h.
 \end{equation}

For finite element discretization in time, we divide the time interval $[0,T]$ into subintervals $I_n=[t_{n-1}, t_n]$, $n = 1,\ldots,N$, with $t_n=n\tau$ 
and stepsize $\tau=T/N $.
For an integer $k \geqslant 1$, we denote by $\mathbb{P}^k$ the space of real-valued polynomials of degree up to $k$ in the time variable $t$. 
Given a Banach space $X$, such as $X = L^2(\Omega)$ or $X = S_h$, the  tensor-product space $\mathbb{P}^k \otimes X$ can be defined as follows:
\[
\mathbb{P}^k \otimes X := \text{span}\big\{p(t)\phi(x) : p \in \mathbb{P}^k, \phi \in X\big\} 
= \Big\{ \sum_{j=0}^{k} t^j \phi_j : \phi_j \in X  \Big\} .
\]

For $n = 1, 2, \ldots, N$ and given $u_h^{n-1}$, we consider the following space-time FEM on $\Omega\times I_n$: Find  
$u_h|_{I_n} \in \mathbb{P}^k\otimes S_h$ and $\kappa_{0,h}|_{I_n},\kappa_{j+1,h}|_{I_n}\in\mathbb{R}$, $j=0,1,\dotsc,d$,  satisfying the following equations:
    \begin{subequations}\label{algorithm}
      \begin{align}
          &\text{Re} \int_{t_{n-1}}^{t_n} \int_{\Omega}\i\partial_t u_h \cdot \i\bar{u}_h \, \d x \d t 
          = 0,  \label{2a}\\
    &\text{Re} \int_{t_{n-1}}^{t_n} \int_{\Omega}\i\nabla\partial_t u_h \cdot \i\nabla\bar{u}_h \, \d x \d t 
    -\text{Re} \int_{t_{n-1}}^{t_n} \int_{\Omega}\i\partial_tu_h\cdot \i f(|u_h|^2)\bar{u}_h  \, \d x \d t
          = 0,  \label{2aa}\\
          &\text{Re} \int_{t_{n-1}}^{t_n} \int_{\Omega}\i\partial_t u_h  \cdot \partial_j \bar{u}_h \, \d x \d t
          = 0\quad\text{for}\,\,\, j = 1,\dotsc,d, \label{2b}\\
& {\rm Re}\int_{t_{n-1}}^{t_n}\int_{\Omega}\i \partial_t u_h \cdot \bar{v}_h\, \d x \d t 
 - {\rm Re}\int_{t_{n-1}}^{t_n} \int_{\Omega}  \nabla u_h \cdot \nabla \bar{v}_h\, \d x \d t 
 + {\rm Re}\int_{t_{n-1}}^{t_n}\int_{\Omega}f(|u_h|^2) u_h \bar{v}_h \, \d x\d t \nonumber \\
= &\, \kappa_{0,h}\,{\rm Re}\int_{t_{n-1}}^{t_n}  \int_{\Omega} \i u_h \bar{v}_h \, \d x \d t + 
\kappa_{1,h}\,{\rm Re}\int_{t_{n-1}}^{t_n}  \int_{\Omega}\i\nabla u_h \nabla\bar{v}_h \, \d x \d t\nonumber \\
&-\kappa_{1,h}\,{\rm Re}\int_{t_{n-1}}^{t_n}  \int_{\Omega}\i f(|u_h|^2)u_h\bar{v}_h \, \d x \d t
\nonumber \\
&+\sum_{j=1}^d \kappa_{j+1,h}\,{\rm Re}\int_{t_{n-1}}^{t_n}\int_{\Omega}  \partial_j u_h \bar{v}_h \, \d x \d t
\quad\forall\, v_h\in \mathbb{P}^{k-1}\otimes S_h , \label{2c} \\
          &u_h(t_{n-1}) = u_h^{n-1} . \hspace{41pt}  \label{2d}
      \end{align}
  \end{subequations}
As initial value of the numerical solution we can take $u^0_h := I_h u_0$, where $I_h$ is the Lagrange interpolation operator onto the finite element space. In \eqref{algorithm}, we view  the coefficients $\kappa_{j,h}$ as piecewise constant functions, constant in each subinterval $I_n.$ 
  
The number of unknown functions $u_h|_{I_n}$, $\kappa_{0,h}|_{I_n}$, and $\kappa_{j+1,h}|_{I_n}$, $j=0,1,\dotsc,d$, is equal to the number of equations 
in \eqref{algorithm}. The constraint equations, namely, \eqref{2a}, \eqref{2aa}, and \eqref{2b}, ensure the conservation of mass, energy, 
and momentum in view of \eqref{mass-conservation}--\eqref{momentum-conservation}. This is presented in the following theorem. 

\begin{theorem}[Conservation properties of the algorithm]\label{th3.1}
  Let $u_h^n\in S_h$, $n=1,\dotsc,N$, be the numerical solutions determined by the algorithm \eqref{algorithm}.  
  Then, the following conservation properties hold:
  \begin{subequations}
    \begin{align}
      M[u_h(t_n)] &= M[u_h(t_{n-1})],\\
            P[u_h(t_n)] &= P[u_h(t_{n-1})] ,\\
      E[u_h(t_n)] &= E[u_h(t_{n-1})] .
    \end{align}
  \end{subequations}
\end{theorem}

\subsection{Equivalent formulation and iterative scheme}
The algorithm \eqref{algorithm} can also be formulated without Lagrange multipliers  by introducing the subspace 
  \begin{align*}
    X_h(u_h|_{I_n}) = {\rm span}\{\i u_h, \i\Delta_h u_h + \i P_h[f(|u_h|^2)u_h], \partial_1u_h,\dotsc, \partial_du_h\} , 
  \end{align*}
 where $P_h: L^2(\Omega)\rightarrow S_h$ stands for the $L^2$-orthogonal projection onto the finite element space. 
 The orthogonal complement of $ X_h(u_h|_{I_n}) $ with respect to the inner product of $\mathbb{P}^{k-1} \otimes S_h$ is defined as 
     \begin{equation*}
    X_h(u_h|_{I_n})^\perp 
    =\Big\{v_h\in \mathbb{P}^{k-1} \otimes S_h: {\rm Re}\int_{t_{n-1}}^{t_n} \int_\Omega v_h\overline\phi_h\ \d x\d t =0\,\,\text{for all}\,\,\phi_h\in X_h(u_h|_{I_n}) \Big\},
     \end{equation*} 
     where the inner product is defined as the real part of the $L^2$ inner product on a complex space, which can be interpreted as an inner product on the corresponding two-dimensional real vector space.

    Then, \eqref{algorithm} can be equivalently reformulated as follows: 
For given $u_h^{n-1}$, find $u_h|_{I_n}\in \mathbb{P}^{k}\otimes S_h$ with $\i\partial_t u_h|_{I_n}\in X_h(u_h|_{I_n})^\perp$ 
and initial condition $ u_h(t_{n-1}) = u_h^{n-1} $, satisfying the following equation: 
\begin{equation}\label{2c-2}  
  \begin{split}
 &{\rm Re}\int_{t_{n-1}}^{t_n}\int_{\Omega}\i \partial_t u_h \cdot \bar{v}_h\, \d x \d t 
 - {\rm Re}\int_{t_{n-1}}^{t_n} \int_{\Omega}  \nabla u_h \cdot \nabla \bar{v}_h\, \d x \d t\\
 & + {\rm Re}\int_{t_{n-1}}^{t_n}\int_{\Omega}f(|u_h|^2) u_h \bar{v}_h \, \d x \d t 
= 0\quad\forall\, v_h\in X_h(u_h|_{I_n})^\perp .
\end{split}
\end{equation}      
We emphasize that the test and trial space $ X_h(u_h|_{I_n})^\perp$ in \eqref{2c-2} depends on the numerical solution $u_h.$ 



For the numerical implementation, the solution of the nonlinearly implicit scheme \eqref{algorithm} needs to be approximated by an iterative algorithm. 
In view of the formulation in \eqref{2c-2} (which is essentially a nonlinear system for $u_h$ without considering $\kappa_{j,h}$), 
we propose the following iterative approach: Newton's iteration is used for $u_h$, while $\kappa_{j,h}$ is treated as an unknown Lagrange multiplier during each iteration of $u_h$, without applying Newton's iteration to it. This is equivalent to applying a fixed-point iteration for the Lagrange multipliers. 
Namely, for given $u_h^{(\ell-1)}|_{I_n}\in\mathbb{P}^k\otimes S_h$, find $u_h^{(\ell)}|_{I_n}\in\mathbb{P}^k\otimes S_h$ and $\kappa_{j,h}^{(\ell)} \in\R$, $j=0,1,\dotsc,d+1$, 
    satisfying the following linearized equations: 
\begin{subequations}\label{4}
      \begin{align}
          &\text{Re}\int_{t_{n-1}}^{t_n}\int_{\Omega}\i\partial_t u_h^{(\ell)} \cdot \i\bar{u}_h^{(\ell-1)}  \, \d x \,\d t+\text{Re}\int_{t_{n-1}}^{t_n}\int_{\Omega}\i\partial_t u_h^{(\ell-1)} \cdot \i(\bar{u}_h^{(\ell)} - \bar{u}_h^{(\ell-1)} ) \, \d x \,\d t = 0, \label{4a}\\
           &\text{Re}\int_{t_{n-1}}^{t_n}\int_{\Omega}\left(\i\nabla\partial_t u_h^{(\ell)} \cdot \i\nabla\bar{u}_h^{(\ell-1)}
           -\i\partial_t u_h^{(\ell)} \cdot \i f(|u_h^{(\ell-1)}|^2)\bar{u}_h^{(\ell-1)}\right)  \, \d x \, \d t \nonumber\\
           & \, +\text{Re}\int_{t_{n-1}}^{t_n}\int_{\Omega}\i\nabla\partial_t u_h^{(\ell-1)} \cdot \i\nabla(\bar{u}_h^{(\ell)}-\bar{u}_h^{(\ell-1)}) \, \d x \, \d t \nonumber\\
           & \, - \text{Re}\int_{t_{n-1}}^{t_n}\int_{\Omega}\i \partial_t u_h^{(\ell-1)} \cdot \i\left(\bar  g_2(u_h^{(\ell-1)})(u_h^{(\ell)}  - u_h^{(\ell-1)}) + g_1(u_h^{(\ell-1)})(\bar{u}_h^{(\ell)} - \bar{u}_h^{(\ell-1)})\right) \, \d x \, \d t = 0,\label{4aaa}\\
          &\text{Re}\int_{t_{n-1}}^{t_n}\int_{\Omega}\i\partial_t u_h^{(\ell)} \cdot \partial_j\bar{u}_h^{(\ell-1)} \, \d x \, \d t +\text{Re}\int_{t_{n-1}}^{t_n}\int_{\Omega}\i\partial_t u_h^{(\ell-1)} \cdot (\partial_j\bar{u}_h^{(\ell)} - \partial_j\bar{u}_h^{(\ell-1)}) \, \d x \, \d t = 0, \label{4b}\\
          &{\rm Re}\int_{t_{n-1}}^{t_n} \int_{\Omega} \i \partial_t u_h^{(\ell)} \cdot \bar{v}_h\,\d x \, \d t 
          - {\rm Re}\int_{t_{n-1}}^{t_n}\int_{\Omega}\nabla u_h^{(\ell)} \cdot \nabla \bar{v}_h\, \d x\, \d t\nonumber\\ 
          &\,+ {\rm Re}\int_{t_{n-1}}^{t_n}\int_{\Omega}f(|u_h^{(\ell-1)}|^2) u_h^{(\ell-1)} \bar{v}_h \, \d x\, \d t 
          + {\rm Re}\int_{t_{n-1}}^{t_n}\int_{\Omega} g_1(u_h^{(\ell-1)})(u_h^{(\ell)} - u_h^{(\ell-1)})\bar{v}_h\, \d x\, \d t \nonumber \\
         &+ {\rm Re}\int_{t_{n-1}}^{t_n}\int_{\Omega}  g_2(u_h^{(\ell-1)})(\bar{u}_h^{(\ell)} - \bar{u}_h^{(\ell-1)}) \bar{v}_h \, \d x\, \d t \nonumber \\
          = &\,\kappa_{0,h}^{(\ell)} \,{\rm Re}\int_{t_{n-1}}^{t_n} \int_{\Omega} \i u_h^{(\ell-1)} \bar{v}_h \, \d x \, \d t 
           +\kappa_{1,h}^{(\ell)} \,{\rm Re}\int_{t_{n-1}}^{t_n} \int_{\Omega} \i\nabla u_h^{(\ell-1)} \nabla \bar{v}_h \, \d x \, \d t\nonumber \\
            &\,-\kappa_{1,h}^{(\ell)} \,{\rm Re}\int_{t_{n-1}}^{t_n} \int_{\Omega} \i f(|u_h^{(\ell-1)}|^2)u_h^{(\ell-1)} \bar{v}_h \, \d x \, \d t \nonumber \\
          &\,+ \sum_{j=1}^d  \kappa_{j+1,h}^{(\ell)} \,{\rm Re}\int_{t_{n-1}}^{t_n}\int_{\Omega}  \partial_j u_h^{(\ell-1)} \bar{v}_h \, \d x \, \d t 
          \quad\forall\, v_h\in \mathbb{P}^{k-1}\otimes S_h , \label{4c}\\
          &u_h^{(\ell)}(t_{n-1}) = u_h^{n-1} , \label{4d}
      \end{align}
\end{subequations}
  where 
    \begin{equation*}
g_1(u) := \partial_u [f(|u|^2)u]\quad\text{and}\quad g_2(u) := \partial_{\bar{u}}[f(|u|^2)u].
    \end{equation*}
The iteration for $\ell$ can be terminated once the error reaches a predefined tolerance level.

\begin{remark}\upshape
Since $u_h|_{I_n}$ is a polynomial of degree $k$ with respect to $t$, the integrals in \eqref{2a} and \eqref{2b} can be evaluated 
exactly with the $k$-point Gauss quadrature, which is exact for polynomials of degree up to $2k-1,$ for discretizing the equations 
in time, while the integrands $\i\partial_t u_h \cdot \i\bar{u}_h$ and $\i\partial_t u_h \cdot \partial_j\bar{u}_h$ in \eqref{2a} and \eqref{2b}, 
respectively, are polynomials of degree $2k-1$ with respect to $t$. 
The integrals in \eqref{2aa} can be evaluated accurately by utilizing an $m$-point Gauss quadrature with a sufficiently large $m$. 
For example, for the cubic NLS equation which corresponds to the function $f(|u_h|^2)=|u_h|^2$, the integrand in \eqref{2aa} 
is a polynomial of degree $4k-1$ with respect to $t$. In this case, the integral can be evaluated exactly by using the $2k$-point Gauss quadrature.
\end{remark}

\begin{remark}\upshape
  This standard Newton iteration for \eqref{2c-2} converges under the condition that the initial guess lies sufficiently close to the solution of \eqref{algorithm}. In particular, employing the numerical solution at the previous time level \(t_{n-1}\) as the initial guess at \(t_n\) provides a suitable approximation, whose error is of order \(O(\tau)\) and is subsequently reduced to order \(O(\tau^{2\ell})\) after \(\ell\) iterations, in accordance with the quadratic convergence of Newton's method. 
  \end{remark}

\section{Numerical results}
\label{sec4}
In this section, several numerical examples are presented to illustrate the high-order accuracy of the numerical scheme 
and the conservation of mass, energy, and momentum. The numerical experiments are performed by using 
the open-source high-performance finite element software NGSolve; see \cite{scho2014}.

\begin{example} [Simulation of one-dimensional Ma-soliton \cite{Ma1979}] \label{Example1} \upshape 
Consider the periodic initial value problem 
\begin{equation}\label{ex4.1-1}
\left.
\begin{alignedat}{2}
\i\partial_t u +\Delta u + 2|u|^2u &= 0 \quad &&\text{ in }[-L,L] \times (0, T] \\
u|_{t=0} &= u_0 \quad &&\text{ in } [-L,L]
\end{alignedat}
\right\},
\end{equation}
which is an approximation of the NLS equation on the real line $\mathbb{R}$. The soliton solution that Ma derived is, after some simplification 
and choice of the space and time origin,
\begin{equation}
    \label{exx4.1}
u(x,t)=q_0\exp{(2\i q_0^2t)}\left[1+\frac{2m_1(m_1\cos{(4m_1m_2q_0^2t)+\i m_2\sin{(4m_1m_2q_0^2t)}})}{m_2\cosh{(2m_1q_0x)}+ \cos{(4m_1m_2q_0^2t)}}\right],
\end{equation}
where $m_2^2=1 + m_1^2$, an unsteady solution with period $\pi/(2m_1m_2q_0^2)$. The proposed Newton iterative method for $u_h$, 
with fixed-point iteration for $\kappa_{j,h}$, is used to solve the nonlinear system. The iteration is set to stop when the $H^1$ error ($\| u_h^{(\ell)}(t_n) - u_h^{(\ell-1)}(t_n) \|_{H^1}$) is below $10^{-9}$. 
The $L^\infty (0,T,H^1)$ error between the numerical solution and the exact solution \eqref{exx4.1} is measured by 
\begin{equation}
\label{exx4.1-2}
H^1 \ \textrm{error} = {\max_{\substack{1\le n\le N \\1\le j\le k}}}\|u_h(\cdot,t_{nj})-P_{p+2}u(\cdot,t_{nj})\|_{H^1},
\end{equation}
where $P_{p+2}$ denotes the $L^2$ projection to the finite element space of degree $p+2$ (two degrees higher than $S_h$), and $t_{nj}$, $j=1,\dotsc,k$, 
denote the Gauss points on $[t_{n-1},t_n]$. 

In this example, we choose $m_1=1$, $q_0=0.5$, and $L=20$. Taking $T=2$, the time discretization errors are presented in Figure \ref{fig1}, 
where we used finite elements of degree $p=3$ with a sufficiently fine spatial mesh $h=1/100$ so that the error 
due to the spatial discretization is negligibly 
small in observing the temporal convergence rates. From Figure \ref{fig1}, we see that the error of the time discretization is $O(\tau^{k+1})$
 in the $L^\infty (0,T,H^1)$-norm.
The spatial discretization errors are shown in Figure \ref{fig2}, where we chose $k = 3$ with a sufficiently small time stepsize $\tau = 1/100$ 
so that the time discretization error is negligibly small compared to the spatial error. The numerical results in Figure \ref{fig2} show that the 
spatial discretization errors are $O(h^{p})$ in the $L^\infty (0,T,H^1)$-norm. 

\begin{figure}[htp!] 
    \centering
    {\includegraphics[trim = .1cm .1cm .1cm .1cm, clip=true,width=0.45\textwidth,height=0.40\textwidth]{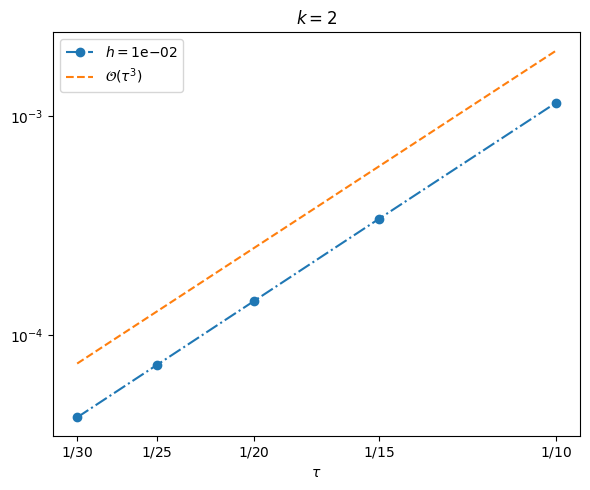}}
    \quad
    {\includegraphics[trim = .1cm .1cm .1cm .1cm, clip=true,width=0.45\textwidth,height=0.40\textwidth]{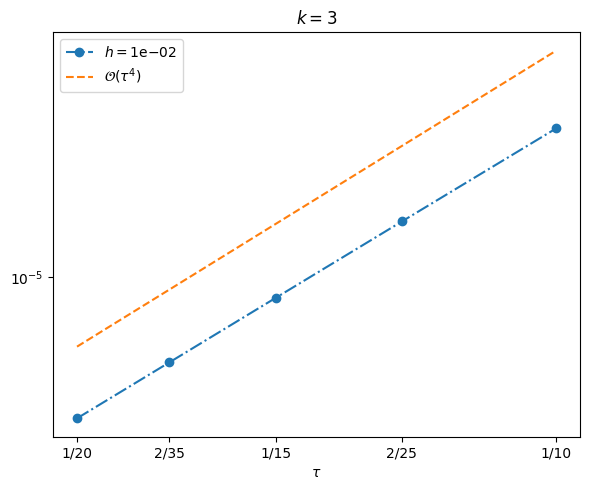}}
    \quad
\vspace{-10pt}
\caption{\small Time discretization errors in the $L^\infty (0,T,H^1)$ norm $($Example \ref{Example1}$)$.}\label{fig1}
\end{figure}

\begin{figure}[htp!] 
    \centering
    {\includegraphics[trim = .1cm .1cm .1cm .1cm, clip=true,width=0.45\textwidth,height=0.40\textwidth]{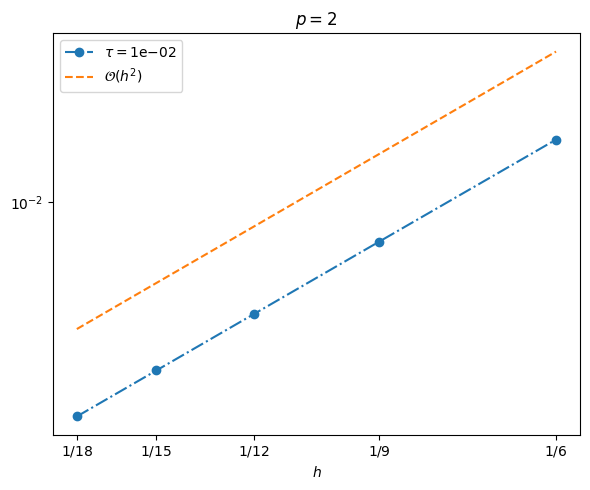}}
    \quad
     {\includegraphics[trim = .1cm .1cm .1cm .1cm, clip=true,width=0.45\textwidth,height=0.40\textwidth]{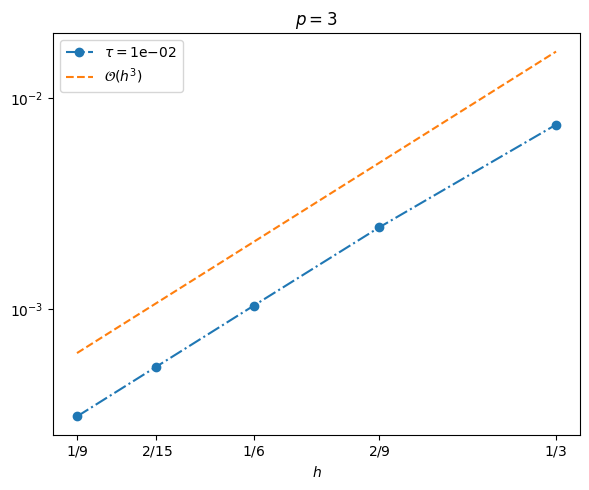}}
\vspace{-10pt}
\caption{\small Space discretization errors in the $L^\infty (0,T,H^1)$ norm $($Example \ref{Example1}$)$.}\label{fig2}
\end{figure}

\begin{figure}[htp!] 
    \centering
    \vspace{10pt}
    {\includegraphics[trim = .1cm .1cm .1cm .1cm, clip=true,width=0.45\textwidth,height=0.40\textwidth]{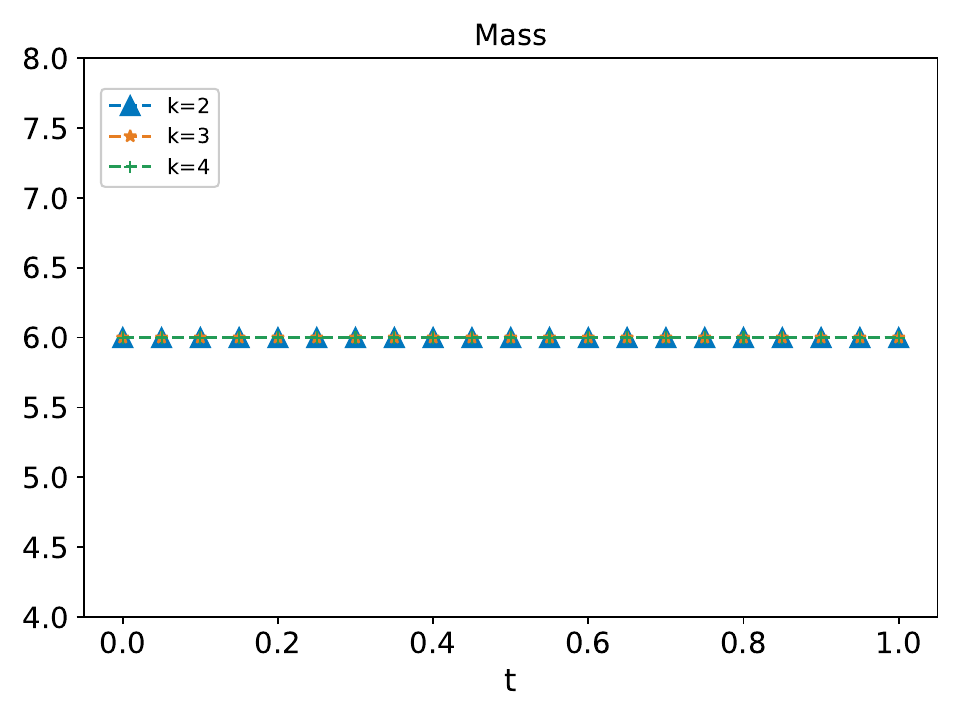}}
    \quad
    {\includegraphics[trim = .1cm .1cm .1cm .1cm, clip=true,width=0.45\textwidth,height=0.40\textwidth]{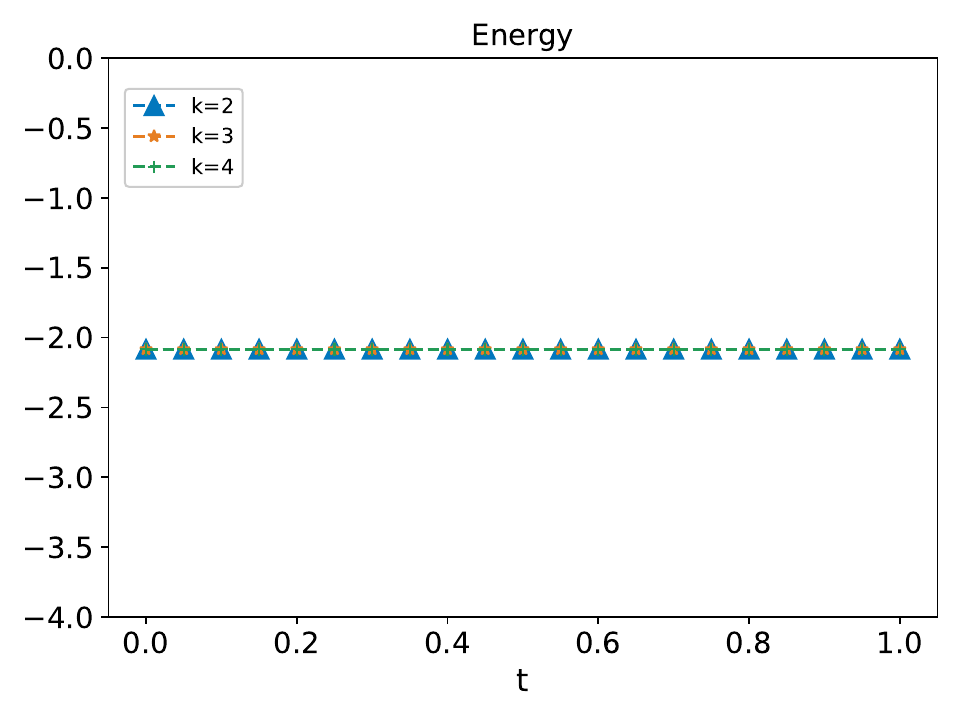}}
    
     {\includegraphics[trim = .1cm .1cm .1cm .1cm, clip=true,width=0.45\textwidth,height=0.40\textwidth]{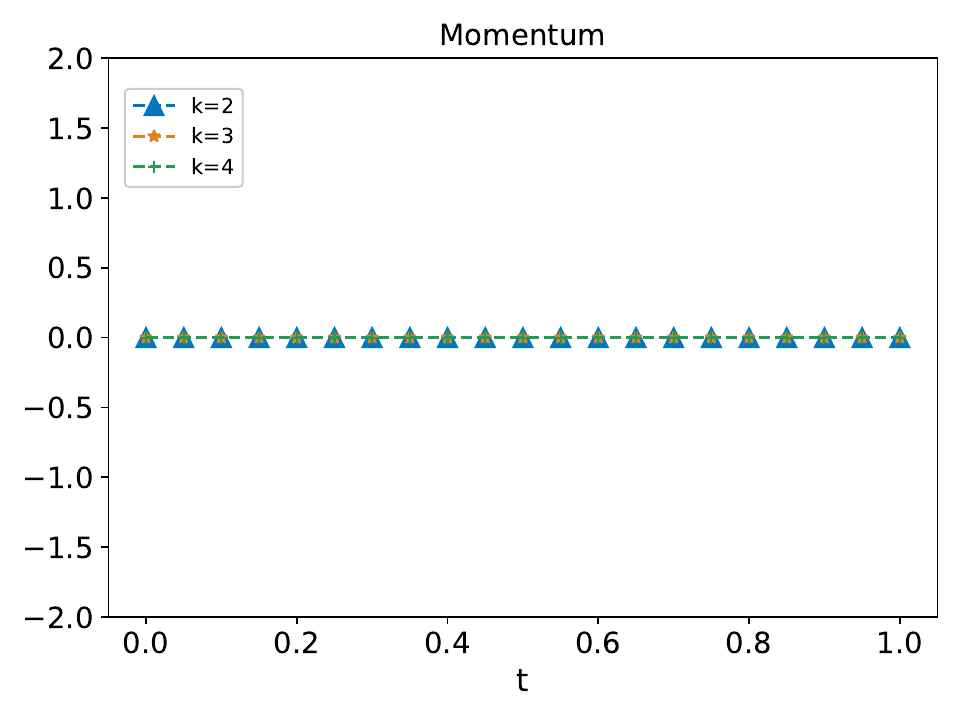}}
\vspace{-10pt}
\caption{\small Evolution of mass, energy, and momentum $($Example \ref{Example1}$)$.}\label{fig4}
\end{figure}

\begin{figure}[htp!] 
    \centering
    \vspace{15pt}
    {\includegraphics[trim = .1cm .1cm .1cm .1cm, clip=true,width=0.41\textwidth,height=0.36\textwidth]{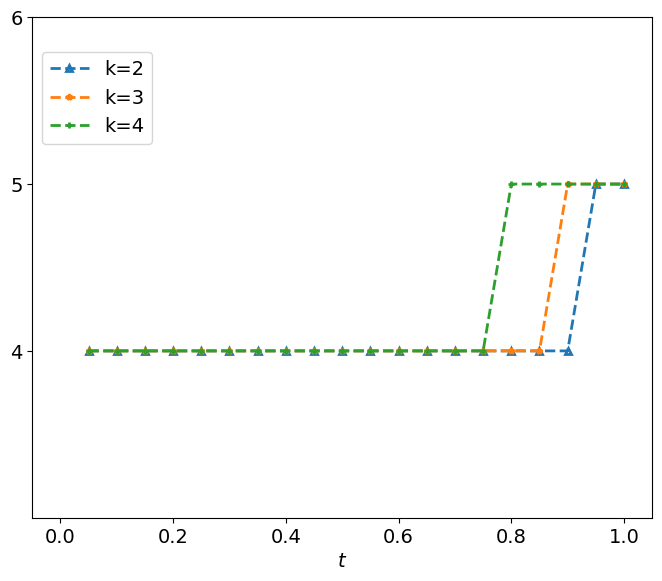}}
     \vspace{-15pt}
   \caption{\small Number of iterations at each time level $($Example \ref{Example1}$)$.}
    \label{fig5}
\end{figure}

\begin{figure}[htp!] 
    \centering
    \vspace{10pt}
    {\includegraphics[trim = .1cm .1cm .1cm .1cm, clip=true,width=0.45\textwidth,height=0.40\textwidth]{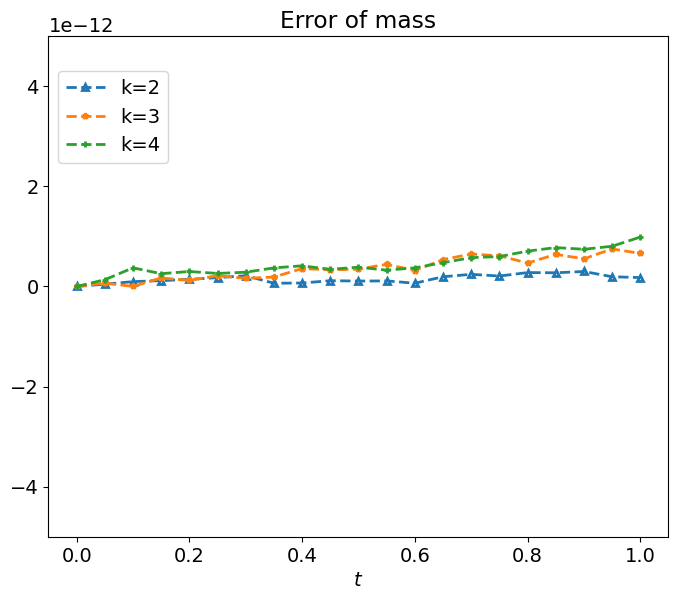}}
    \quad
    {\includegraphics[trim = .1cm .1cm .1cm .1cm, clip=true,width=0.45\textwidth,height=0.40\textwidth]{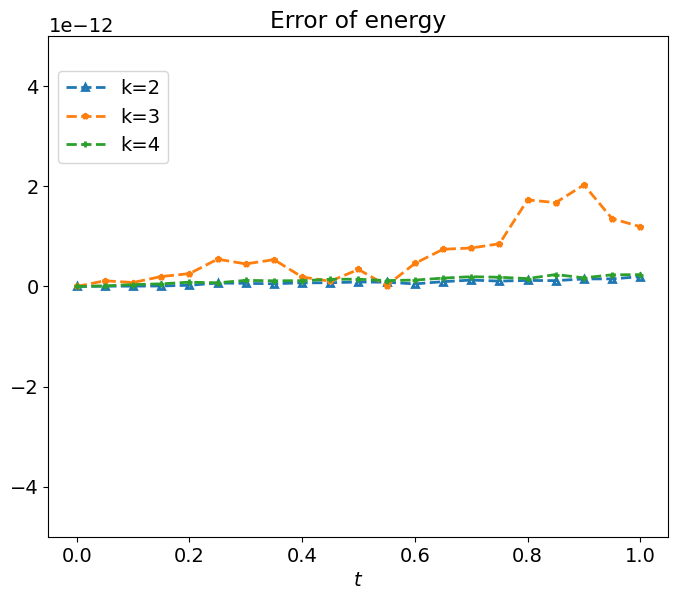}}
    
     {\includegraphics[trim = .1cm .1cm .1cm .1cm, clip=true,width=0.47\textwidth,height=0.40\textwidth]{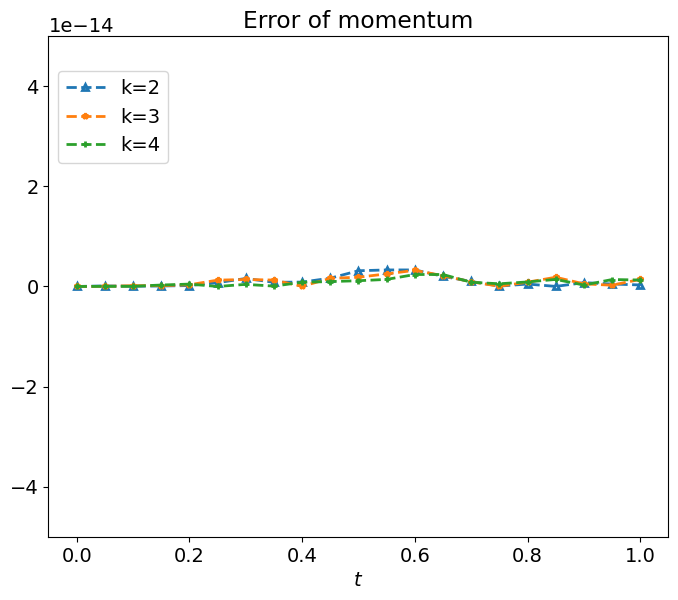}}
     \vspace{-15pt}
\caption{\small Errors of mass, energy, and momentum $($Example \ref{Example1}$)$.}\label{fig3}
\end{figure}

\begin{figure}[htp!] 
  \centering
  {\includegraphics[trim = .1cm .1cm .1cm .1cm, clip=true,width=0.45\textwidth,height=0.40\textwidth]{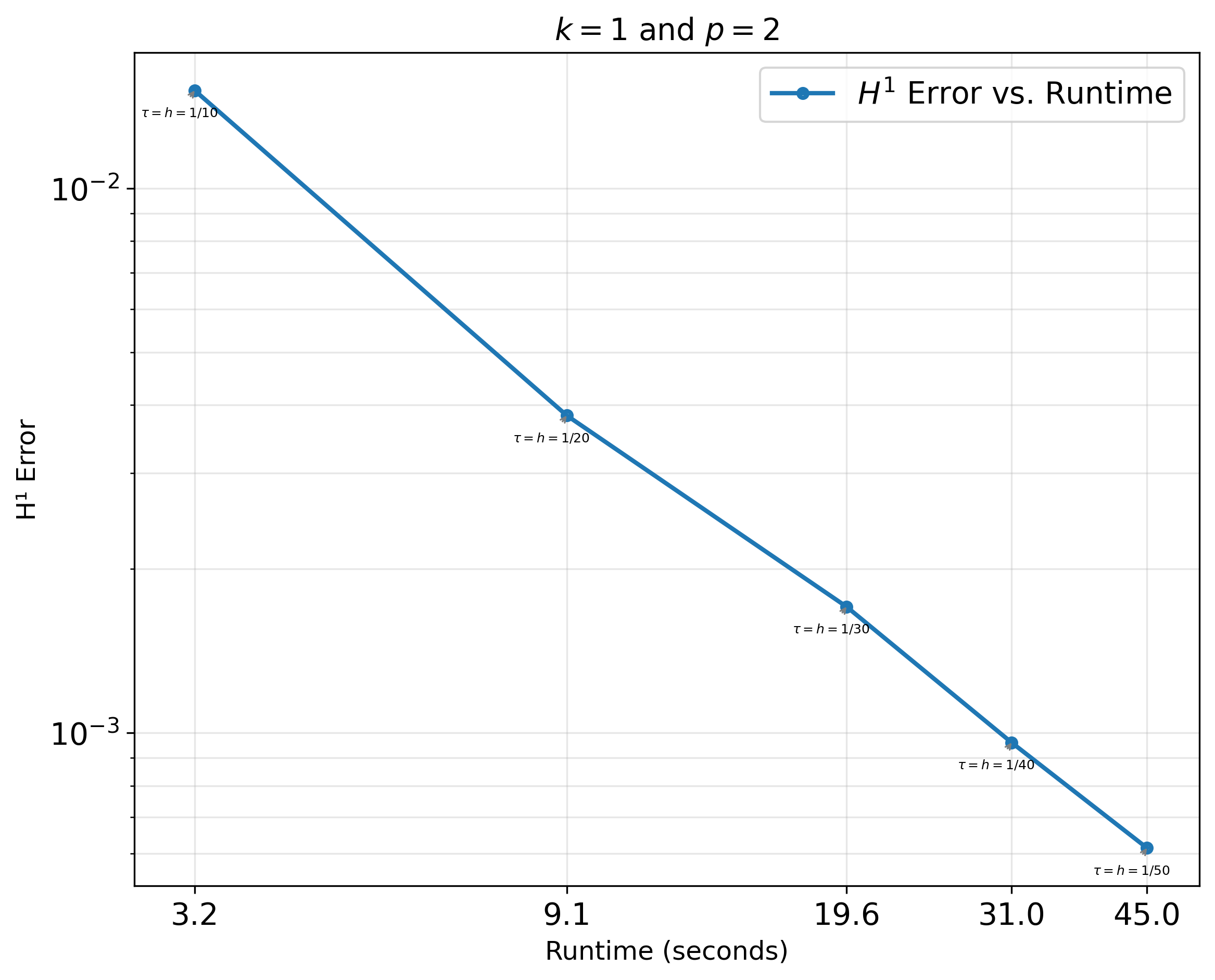}}
  \quad
   {\includegraphics[trim = .1cm .1cm .1cm .1cm, clip=true,width=0.47\textwidth,height=0.40\textwidth]{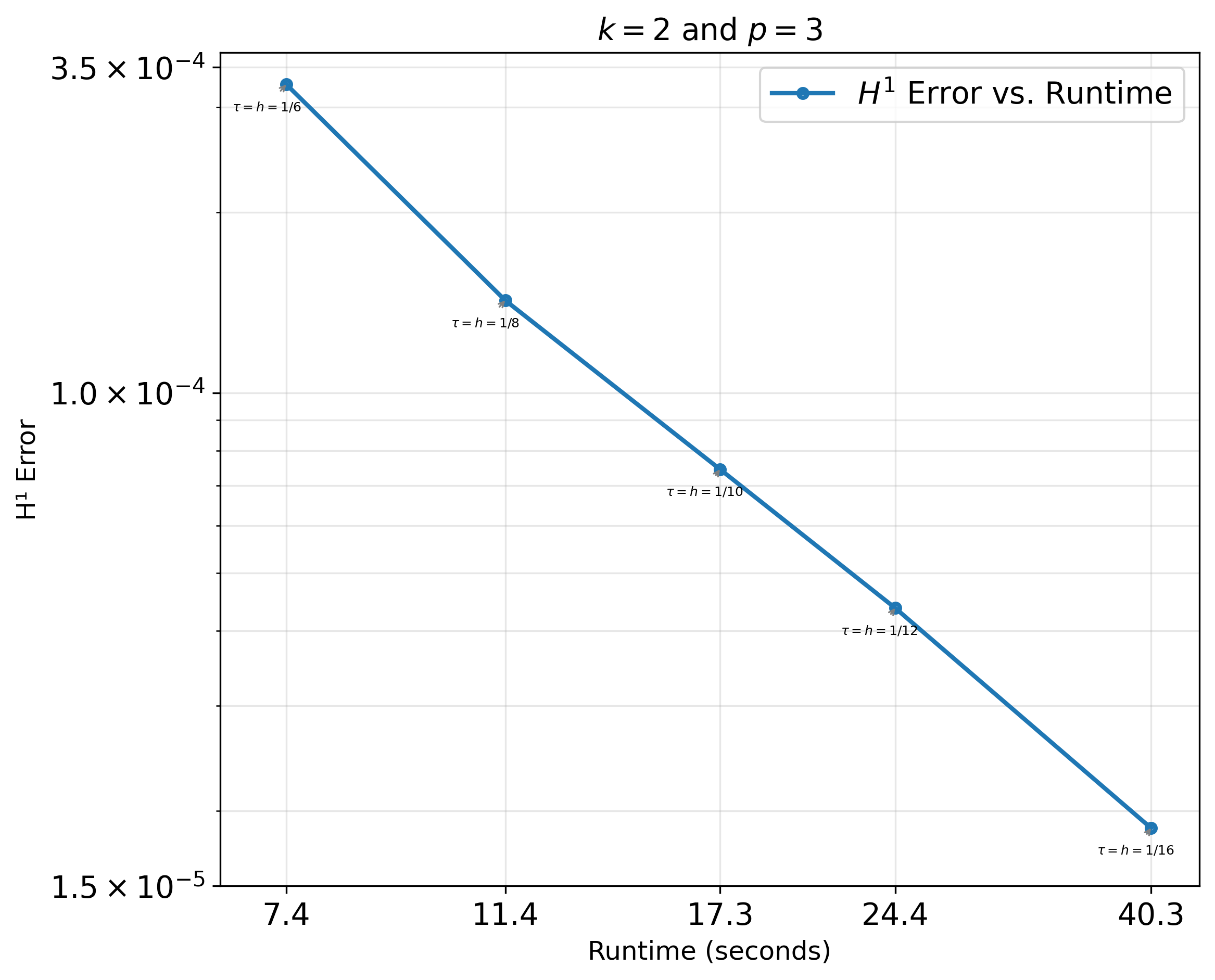}}
\vspace{-10pt}
\caption{\small Change in $H^1$ error as total runtime increases $($Example \ref{Example1}$)$.}\label{fig6}
\end{figure}

Next, for the numerical solution with $p = 3$,  $k = 3$, $\tau=1/20$, and $h = 1/16$, we present the evolutions of mass, energy, and momentum in Figure \ref{fig4}, which shows that there is no visible loss of mass, energy, and momentum in the evolution. 
The errors of mass, energy, and momentum in the numerical solution are presented in Figure \ref{fig3}; we see that the mass, energy, and momentum are conserved up to a discrepancy of the order $10^{-12}$, which is significantly 
smaller than the errors of the numerical solution (about $10^{-4}$ according to the case $p=3$ in Figure \ref{fig2}). Therefore, 
the numerical results show the high-order convergence of the method as well as its effectiveness in conserving mass, energy, 
and momentum. The numbers of iterations at each time level are presented in Figure \ref{fig5}; we see that a few iterations, 
5 or 6, suffice for the conservation of the mass, energy, and momentum of the numerical solution. 

Furthermore, Figure \ref{fig6} illustrates the variation of the $H^1$ error versus total runtime up to time $T=1$, where we have set $k-1=p$ so that the temporal and spatial discretizations converge at the same rate in the $L^\infty(0,T;H^1)$ norm. The results indicate that, as the time step and mesh size decrease, the $H^1$ error decreases as the runtime increases. Moreover, higher-order methods have significantly smaller errors for the same runtime, underscoring the advantages and importance of designing higher-order methods for the NLS equation. This aligns with the main objective of this paper: developing high-order structure-preserving methods. 


\end{example}
\bigskip

\begin{example}[Simulation of one-dimensional bi-soliton \cite{Zakharov-Shabat-1972}] \label{Example2} \upshape 
We consider the one-dimensional focusing NLS equation on the real line $\mathbb{R}$ with the following solution:
\begin{align}
    \label{ex4.1}
u(x,t)=\frac{\ee^{\i M_1^2t}M_1 \textrm{sech}{M_1x}-\ee^{\i M_2^2t}M_2 \textrm{sech}{M_2x}  }{\cosh{J}-\sinh{J}(\tanh{M_1x}\tanh{M_2x} +\cos{S}\textrm{sech}{M_1x}\textrm{sech}{M_2x})},
\end{align}
with
\begin{align*}
S=(M_1^2-M_2^2)t,\quad
\tanh{J}
=2M_1M_2/(M_1^2+M_2^2).
\end{align*}
The solution in \eqref{ex4.1} represents the interaction between two individual solitons; see \cite{Peregrine}. In this example, we choose $L=20$, $M_1= 1.2$, and $M_2 = 1$. 

\begin{figure}[htp] 
    \centering
    \vspace{10pt}
    \subfigure[Initial function and short time evolution]{\includegraphics[trim = .1cm .1cm .1cm .1cm, clip=true,width=0.48\textwidth,height=0.40\textwidth]{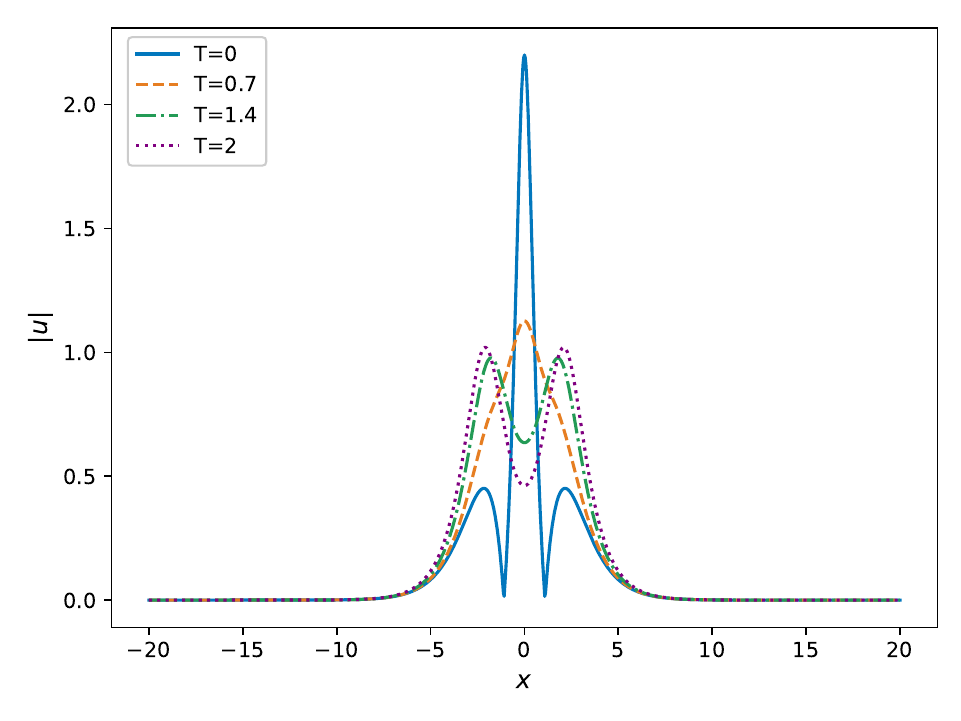}}
    \subfigure[Long time evolution of the numerical solution of the proposed method]{\includegraphics[trim = .1cm .1cm .1cm .1cm, clip=true,width=0.48\textwidth,height=0.40\textwidth]{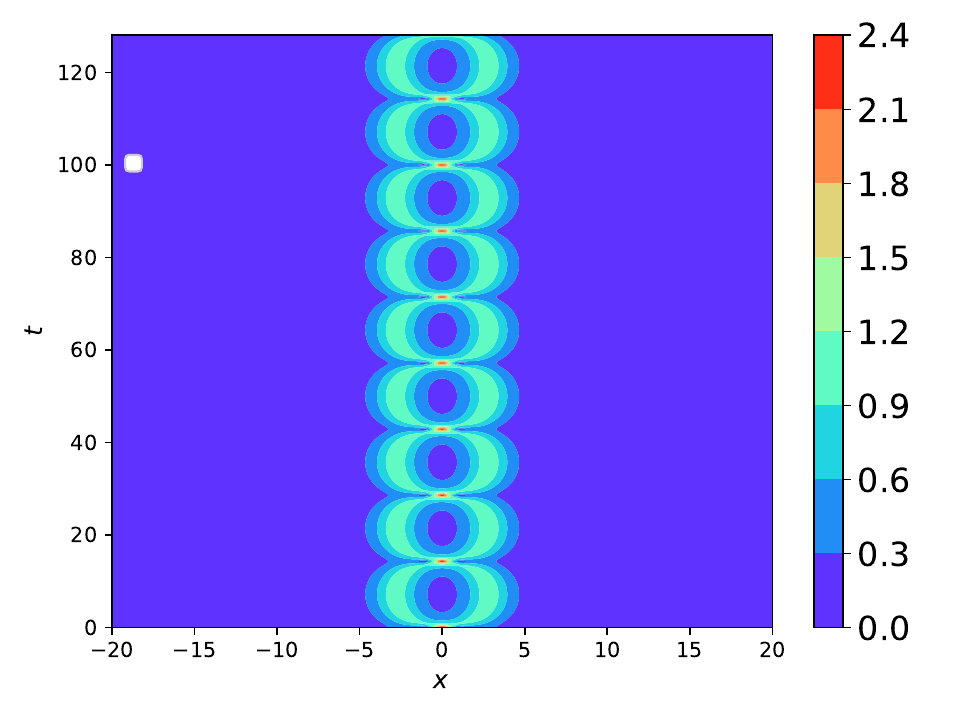}}
    \vspace{-10pt}
\caption{Evolution of the amplitude of the numerical solution $($Example \ref{Example2}$)$.}\label{fig2-6}
\end{figure}

At $t = 0$, the strong interaction between two solitons results in a striking peak of $|u|$ at the origin, as displayed in Figure \ref{fig2-6}(a). 
This peak shows larger $L^\infty$ and $H^1$ norms of the initial function. The most distinctive characteristic of $|u|$ is its periodicity, with 
a period of $2\pi/(M_1^2-M_2^2)$. Consequently, the periodic appearance of the initial peak (see Figure \ref{fig2-6}(b)) presents challenges 
for numerical methods in terms of energy conservation and accuracy, making this situation an ideal choice for examining the long-term 
performance of the proposed method.

\begin{figure}[htp] 
    \centering
    \vspace{10pt}
    {\includegraphics[trim = .1cm .1cm .1cm .1cm, clip=true,width=0.45\textwidth,height=0.40\textwidth]{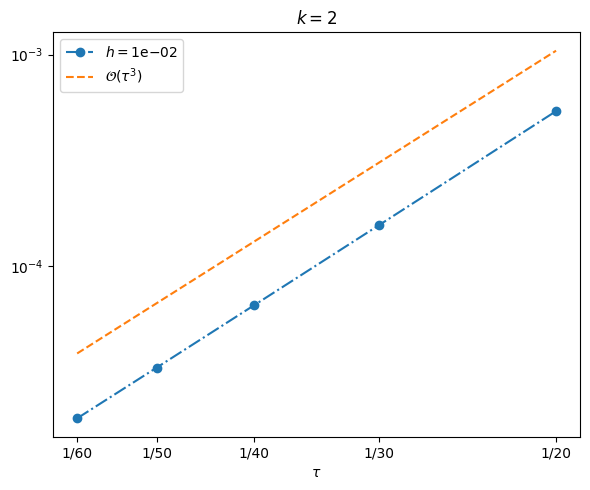}}
    \quad
    {\includegraphics[trim = .1cm .1cm .1cm .1cm, clip=true,width=0.45\textwidth,height=0.40\textwidth]{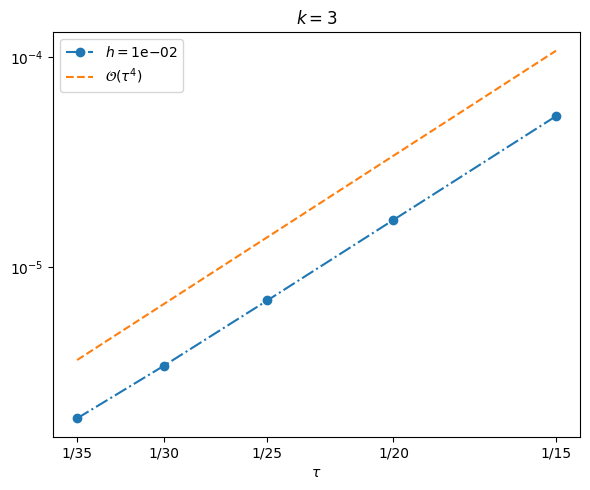}}
\vspace{-10pt}
\caption{\small Time discretization errors in the $L^\infty (0,T,H^1)$ norm $($Example \ref{Example2}$)$.}\label{fig2-1}
\end{figure}
  \begin{figure}[htp!] 
    \centering
    \vspace{10pt}
    {\includegraphics[trim = .1cm .1cm .1cm .1cm, clip=true,width=0.45\textwidth,height=0.40\textwidth]{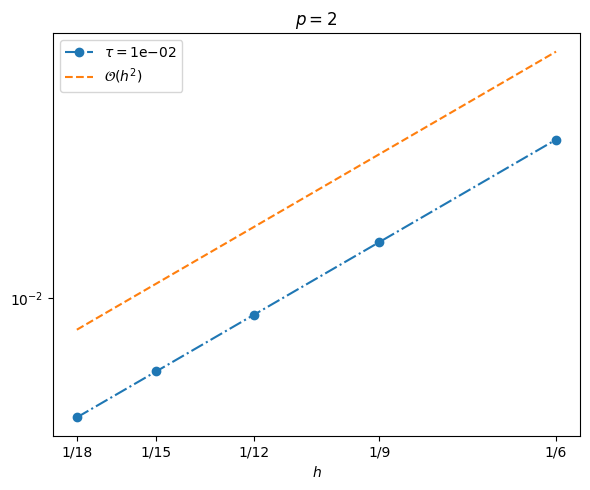}}
    \quad
     {\includegraphics[trim = .1cm .1cm .1cm .1cm, clip=true,width=0.45\textwidth,height=0.40\textwidth]{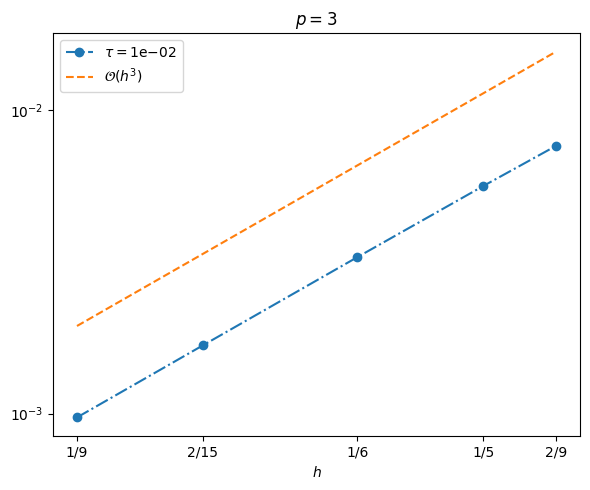}}
\vspace{-10pt}
\caption{\small Space discretization errors in the $L^\infty (0,T,H^1)$ norm $($Example \ref{Example2}$)$.}\label{fig2-2}
\end{figure}

\begin{figure}[htp!] 
    \centering
    \vspace{10pt}
    {\includegraphics[trim = .1cm .1cm .1cm .1cm, clip=true,width=0.45\textwidth,height=0.40\textwidth]{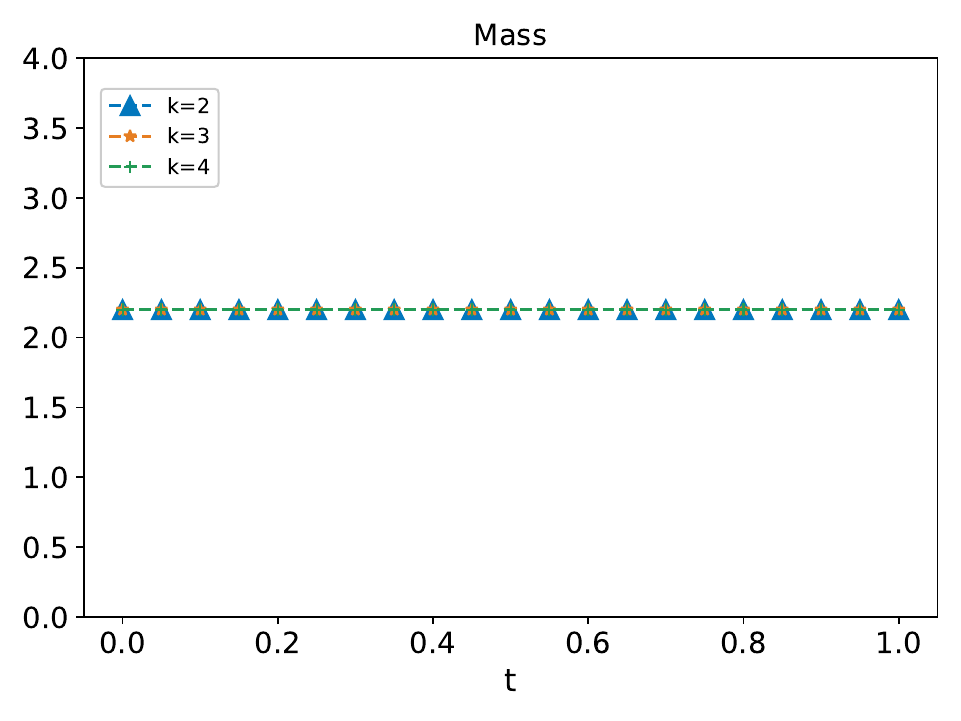}}
    {\includegraphics[trim = .1cm .1cm .1cm .1cm, clip=true,width=0.45\textwidth,height=0.40\textwidth]{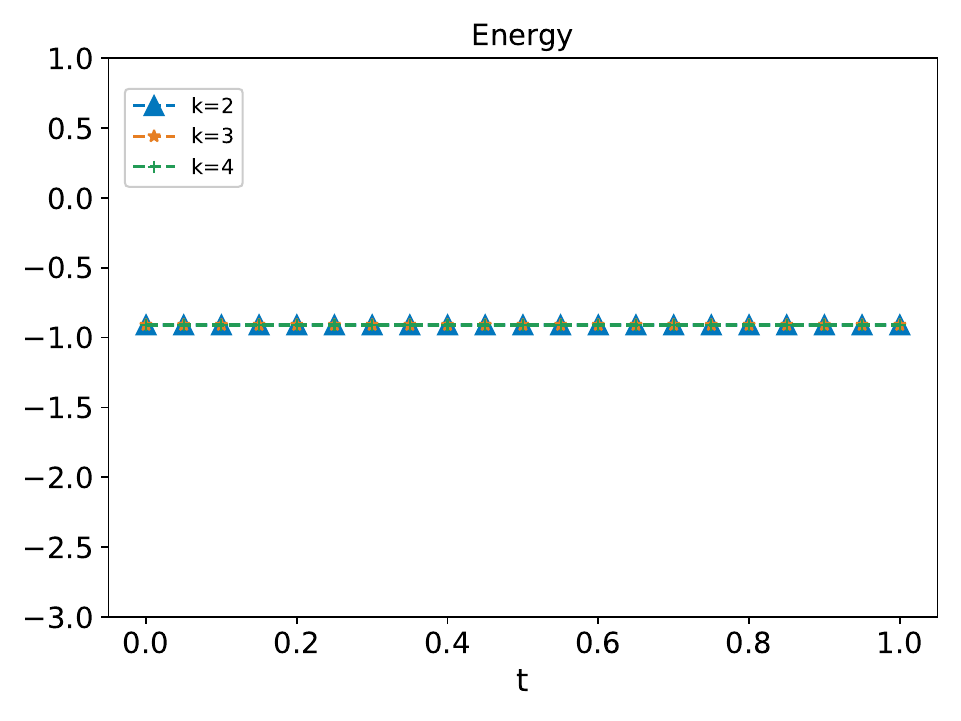}}
    
     {\includegraphics[trim = .1cm .1cm .1cm .1cm, clip=true,width=0.45\textwidth,height=0.40\textwidth]{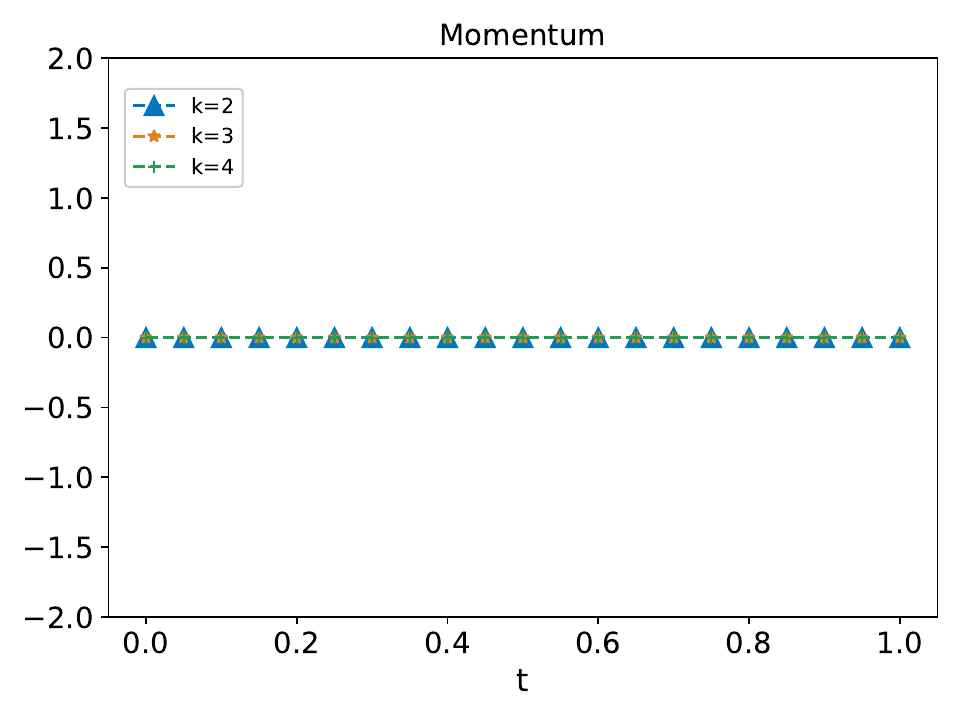}}
\vspace{-10pt}
\caption{\small Time evolutions of mass, energy and  momentum $($Example \ref{Example2}$)$.}\label{fig2-4}
\end{figure}

\begin{figure}[htp!] 
    \centering
    \vspace{10pt}
    {\includegraphics[trim = .1cm .1cm .1cm .1cm, clip=true,width=0.45\textwidth,height=0.4\textwidth]{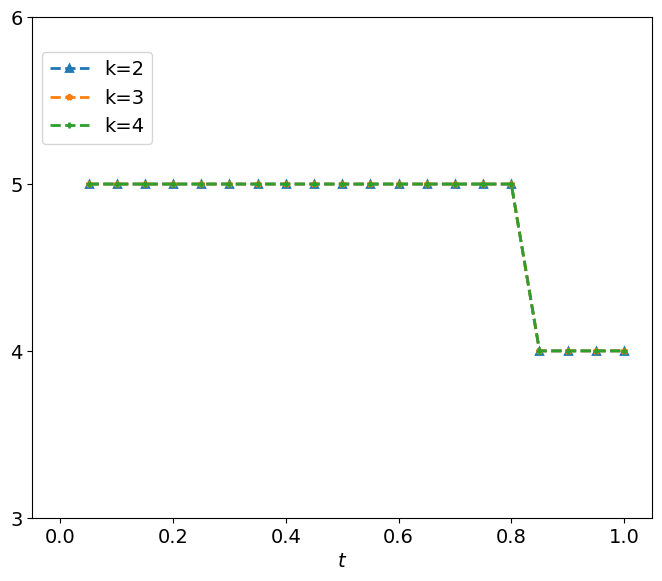}}
\vspace{-10pt}
   \caption{\small Number of iterations at each time level $($Example \ref{Example2}$)$.}
   \label{fig2-5}
\end{figure}

   \begin{figure}[htp!] 
    \centering
    {\includegraphics[trim = .1cm .1cm .1cm .1cm, clip=true,width=0.45\textwidth,height=0.40\textwidth]{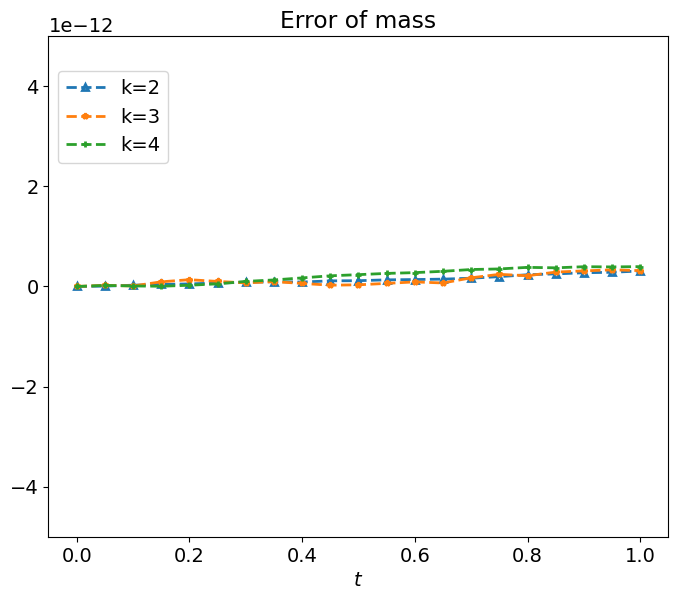}}
    {\includegraphics[trim = .1cm .1cm .1cm .1cm, clip=true,width=0.45\textwidth,height=0.40\textwidth]{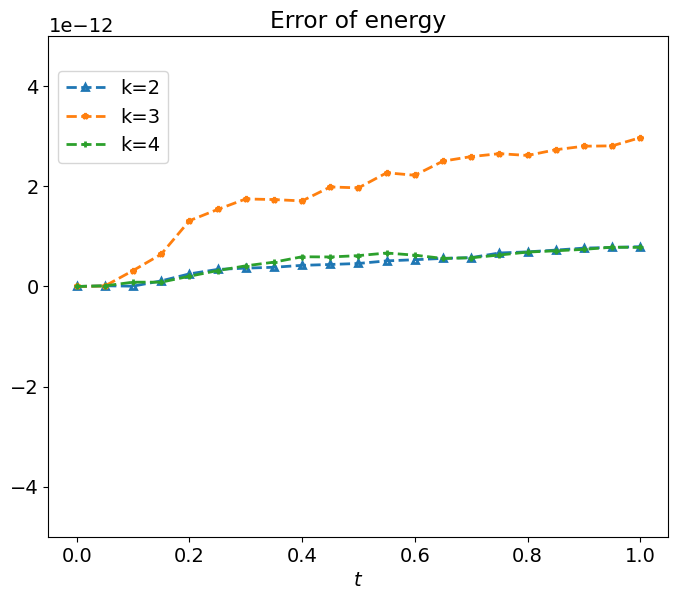}}
     {\includegraphics[trim = .1cm .1cm .1cm .1cm, clip=true,width=0.45\textwidth,height=0.40\textwidth]{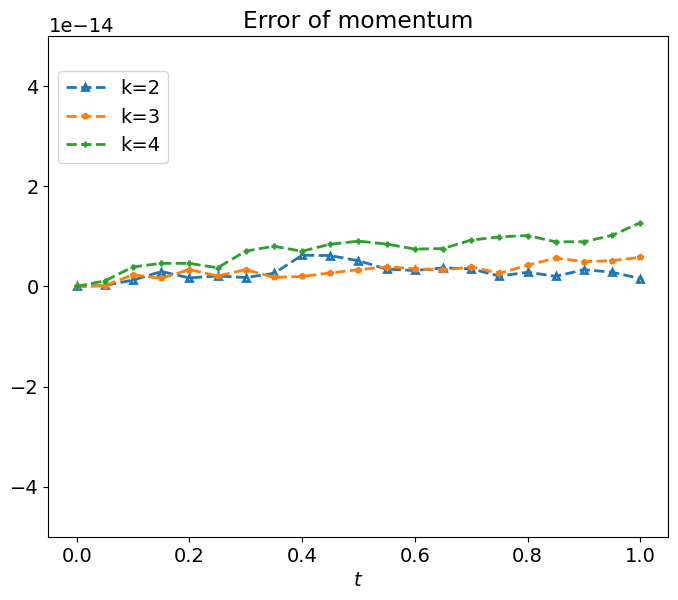}}
\caption{\small Errors of mass, energy, and momentum $($Example \ref{Example2}$)$.}\label{fig2-3}
\end{figure}

\begin{figure}[htp!] 
    \centering
    {\includegraphics[trim = .1cm .1cm .1cm .4cm, clip=true,width=0.75\textwidth,height=0.26\textwidth]{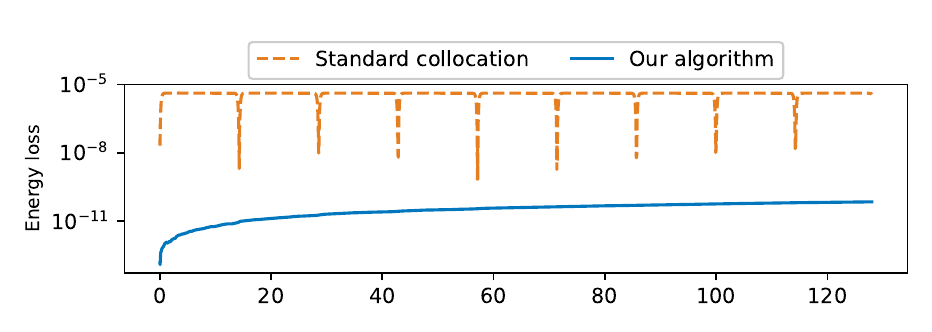}}
    {\includegraphics[trim = .1cm .1cm .1cm .4cm, clip=true,width=0.75\textwidth,height=0.26\textwidth]{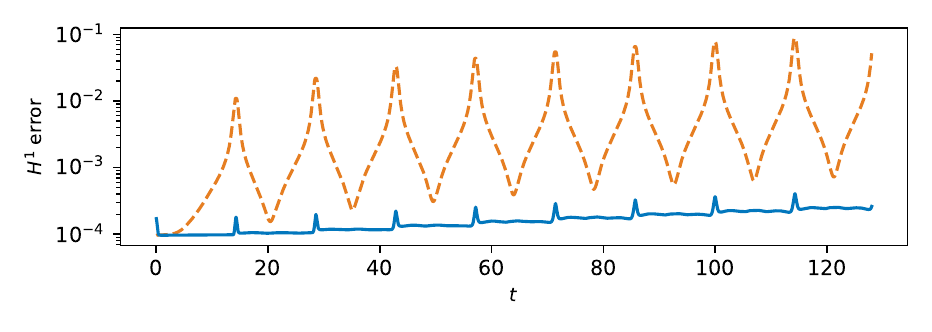}}
\vspace{-10pt}
\caption{\small Energy and $H^1$ errors in long-time simulation, up to $T =128$ $($Example \ref{Example2}$)$.}\label{fig2-66}
\end{figure}

We solve problem \eqref{ex4.1-1} by the proposed method \eqref{4a}--\eqref{4d} and compare the numerical solutions 
with the exact solution \eqref{ex4.1}. The errors at $T=1$ from the temporal discretizations are investigated in Figure \ref{fig2-1} with $p=3$ 
and $h=1/100$ so that the errors caused by the spatial discretization are negligible in observing the temporal convergence orders. 
By fixing $k=3$ and $\tau=1/100$ to neglect the temporal errors in observing the spatial convergence orders, the errors from 
the spatial discretizations are checked in Figure \ref{fig2-2}. The above numerical results indicate that the errors from the temporal 
and spatial discretizations are of order $O (\tau^{k+1})$ and $O (h^p)$,  respectively, in the $L^\infty(0,T,H^1)$-norm.

The time evolution of mass, energy, and momentum is shown in Figure \ref{fig2-4} with $p = 3$, $k = 3$, $T = 1$, 
$\tau = 1/20$, and $h = 1/16$. Figure \ref{fig2-3} shows that the mass, energy, and momentum are conserved 
up to a discrepancy of order $10^{-12}$, which is significantly smaller than the error of the numerical solution 
(about $10^{-7}$ according to Figure \ref{fig2-1}). This shows 
the effectiveness of the proposed method in conserving the mass, energy, and momentum of the NLS equation. 
The numbers of iterations at each time level are shown in Figure \ref{fig2-5}; we see that 6 to 8 iterations suffice 
for  the conservation of the mass, energy, and momentum of the numerical solutions. 

The parameters for the numerical discretization are chosen as $\tau=2^{-5}$, $h=2^{-4}$, $k=2$, and $p=3$, 
with an end time  set to $T=128$. The exact solution in \eqref{ex4.1} decays exponentially as $|x|\rightarrow \infty$, 
and the selected parameter settings guarantee that the solutions before $T=128$ have negligible amplitude 
(up to rounding error) at the boundary of the truncated domain $[-20,20]$. We compared the performance 
of our  algorithm with the standard Gauss collocation method in terms of energy loss and $H^1$ error. 
As shown in Figure \ref{fig2-66}, our algorithm maintains energy conservation within the rounding error range 
and is significantly superior to the standard Gauss collocation method. In addition, our algorithm substantially 
reduces the $H^1$ error of the numerical solutions due to its structure-preserving properties.
\end{example}
\bigskip


\begin{example}[Simulation of two-dimensional soliton]\label{Example3} \upshape 
We investigate the interactions of two-dimensional solitons described by the cubic focusing NLS equation with the following initial value: 
\begin{equation}
\label{ex4.3-2}
u_0(x,y) = \sum_{k=0}^1\exp{(-(x-(-1)^{k}2)^2-y^2)}
\exp{(\i 0.1(x-(-1)^{k}2)+y)}.
\end{equation}
The solution of this problem decays exponentially as $|x + y|$ tends to $\infty$ and is a constant along $x + y = C$. 
Correspondingly, we solve the equation by the proposed method in a rectangular domain (which is periodic in $x\pm y$)
\[\Omega = \{ (x, y) \in \mathbb{R}^2: |x| + |y| \leqslant \sqrt{2}L\} \]
with $L = 10$. Since the exact solution of Example \ref{Example3} is not explicitly given, we compute a reference solution 
$u_{\text{ref}}$ by using a sufficiently small time stepsize $\tau$ for a fixed mesh size $h$ when we test the convergence rates 
of the spatial discretizations, or a sufficiently small mesh size $h$ for a fixed time stepsize $\tau$ when we test the convergence 
rates of the temporal discretizations. 

\begin{figure}[htp!] 
    \centering
    \vspace{-10pt}
    {\includegraphics[trim = .1cm .1cm .1cm .1cm, clip=true,width=0.45\textwidth,height=0.36\textwidth]{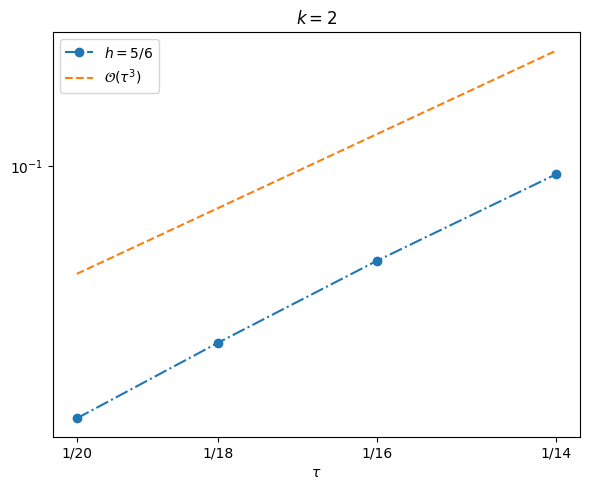}}
    \quad
    {\includegraphics[trim = .1cm .1cm .1cm .1cm, clip=true,width=0.45\textwidth,height=0.36\textwidth]{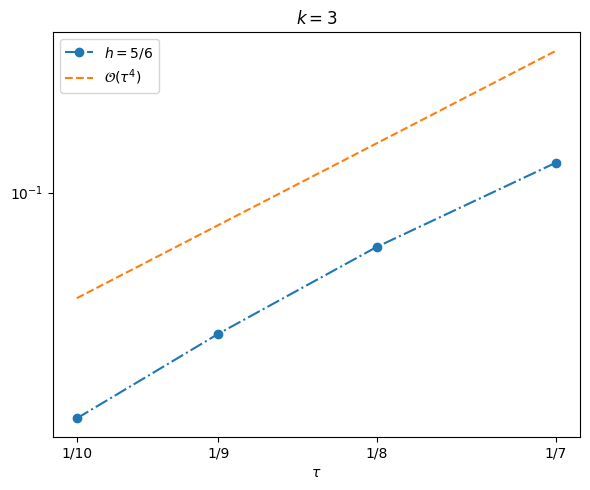}}
\vspace{-10pt}
\caption{\small Time discretization errors in the $L^\infty (0,T,H^1)$  norm $($Example \ref{Example3}$)$.}\label{fig3-1}
\end{figure}
\begin{figure}[htp!] 

\vspace{15pt}
    \centering
    {\includegraphics[trim = .1cm .1cm .1cm .1cm, clip=true,width=0.45\textwidth,height=0.36\textwidth]{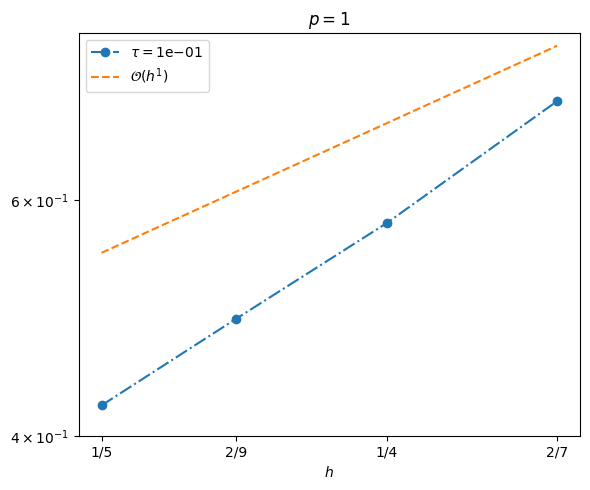}}
    \quad
     {\includegraphics[trim = .1cm .1cm .1cm .1cm, clip=true,width=0.45\textwidth,height=0.36\textwidth]{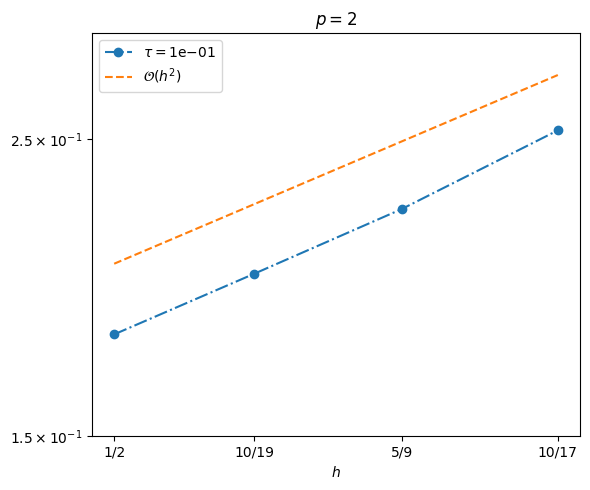}}
\vspace{-10pt}
\caption{\small Space discretization errors in the $L^\infty (0,T,H^1)$ norm $($Example \ref{Example3}$)$.}\label{fig3-2}
\end{figure}

\begin{figure}[htp!]


\vspace{20pt}
    \centering
    {\includegraphics[trim = .1cm .1cm .1cm .1cm, clip=true,width=0.44\textwidth,height=0.36\textwidth]{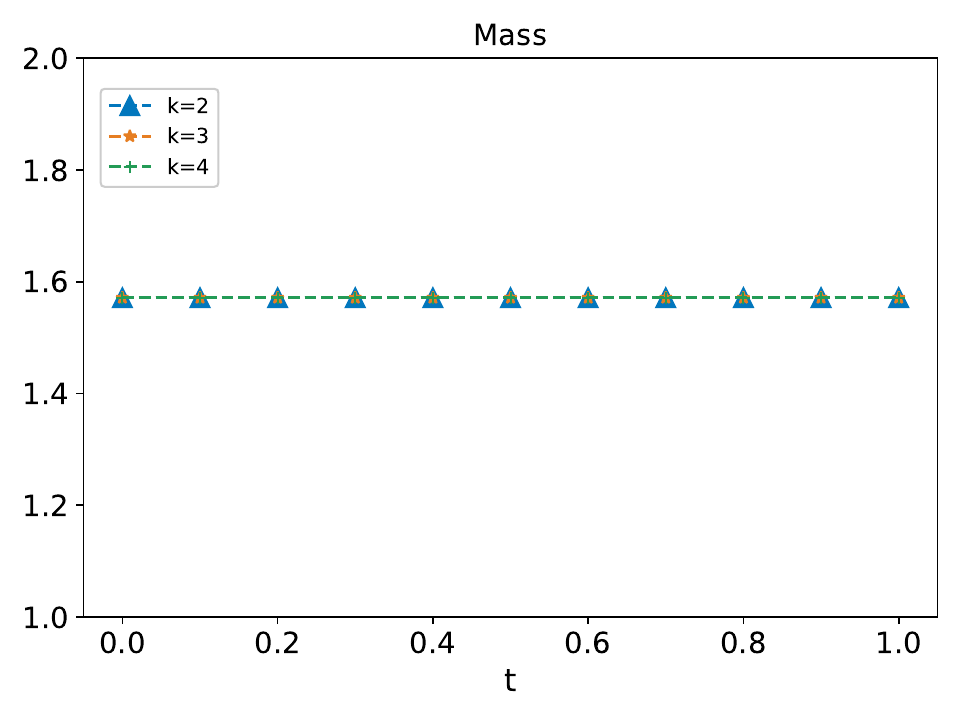}}
    {\includegraphics[trim = .1cm .1cm .1cm .1cm, clip=true,width=0.44\textwidth,height=0.36\textwidth]{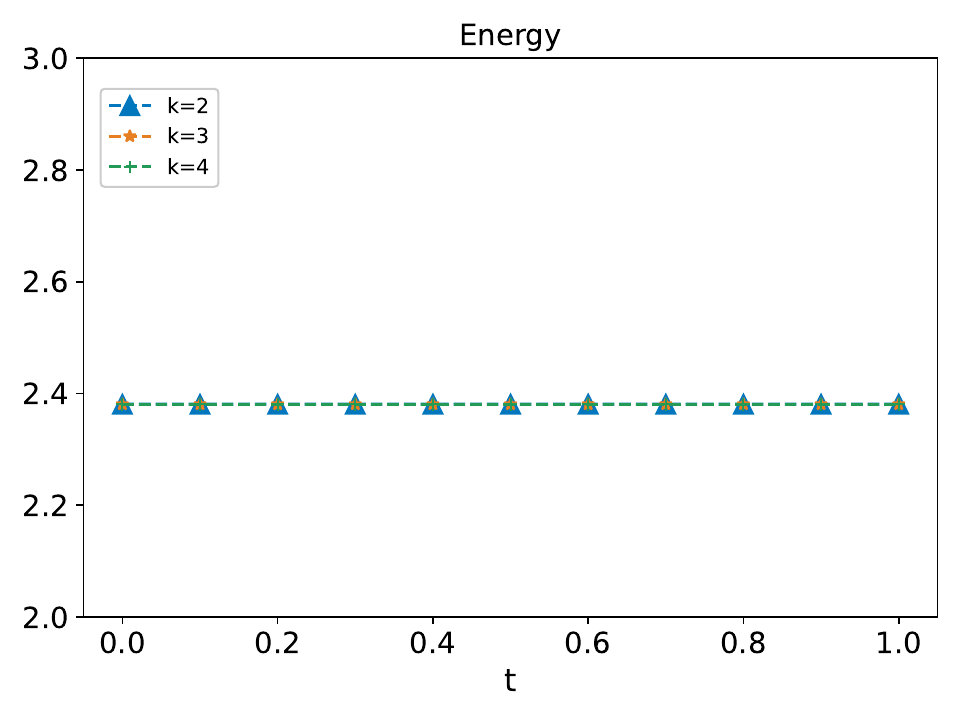}}
     {\includegraphics[trim = .1cm .1cm .1cm .1cm, clip=true,width=0.44\textwidth,height=0.36\textwidth]{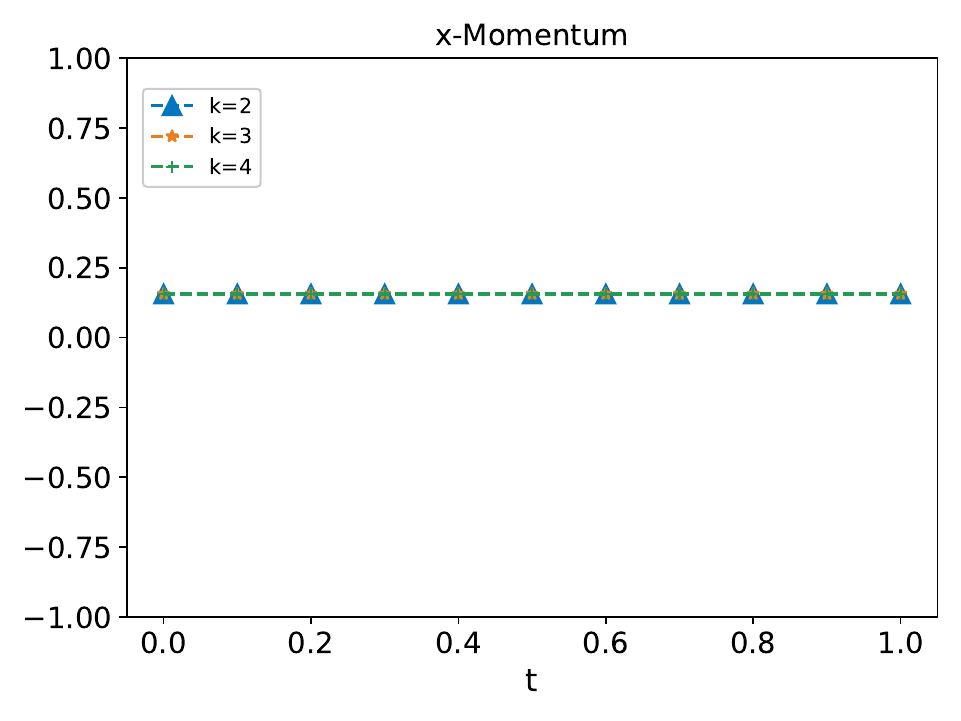}}
     {\includegraphics[trim = .1cm .1cm .1cm .1cm, clip=true,width=0.44\textwidth,height=0.36\textwidth]{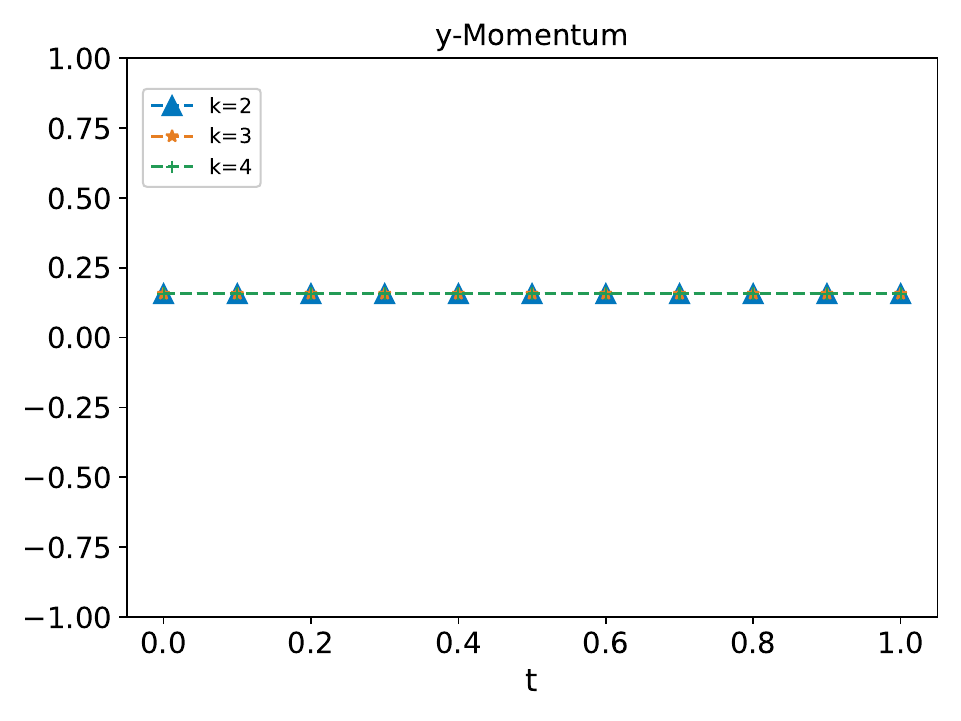}}
\vspace{-10pt}
\caption{\small Time evolutions of mass, energy, and momentum $($Example \ref{Example3}$)$.}\label{fig3-3}

\end{figure}

\begin{figure}[htp!] 
    \centering
    {\includegraphics[trim = .1cm .1cm .1cm .1cm, clip=true,width=0.44\textwidth,height=0.36\textwidth]{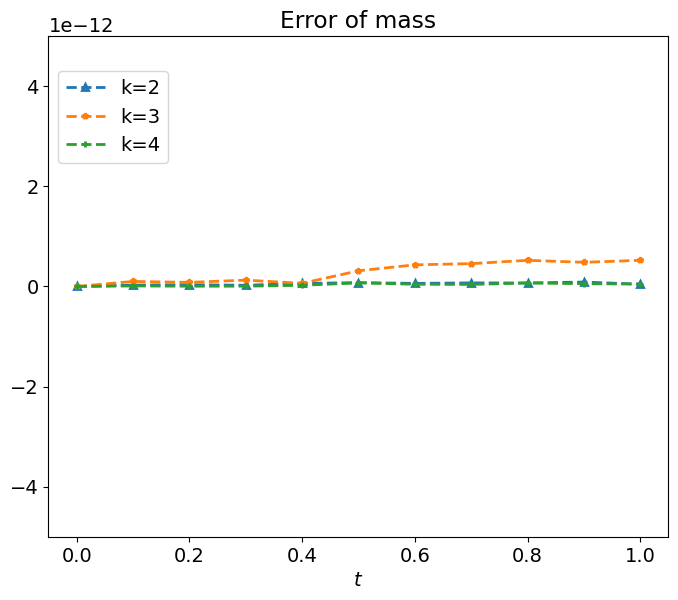}}
    {\includegraphics[trim = .1cm .1cm .1cm .1cm, clip=true,width=0.44\textwidth,height=0.36\textwidth]{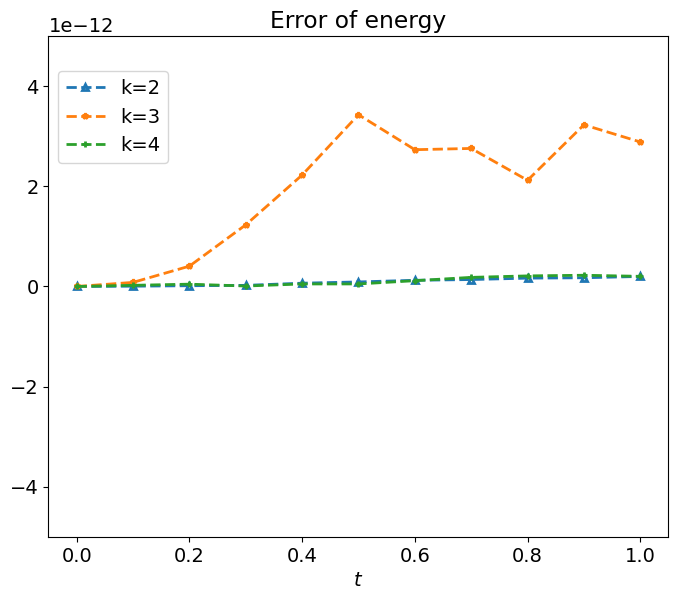}}
     {\includegraphics[trim = .1cm .1cm .1cm .1cm, clip=true,width=0.44\textwidth,height=0.36\textwidth]{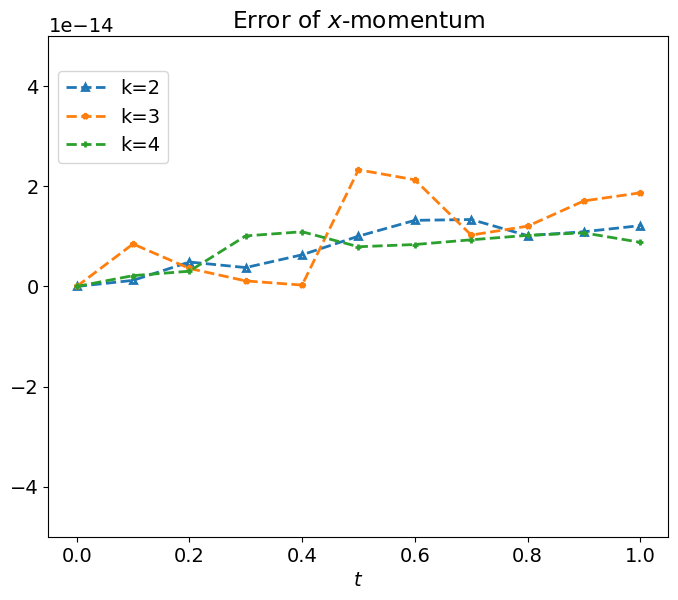}}
      {\includegraphics[trim = .1cm .1cm .1cm .1cm, clip=true,width=0.44\textwidth,height=0.36\textwidth]{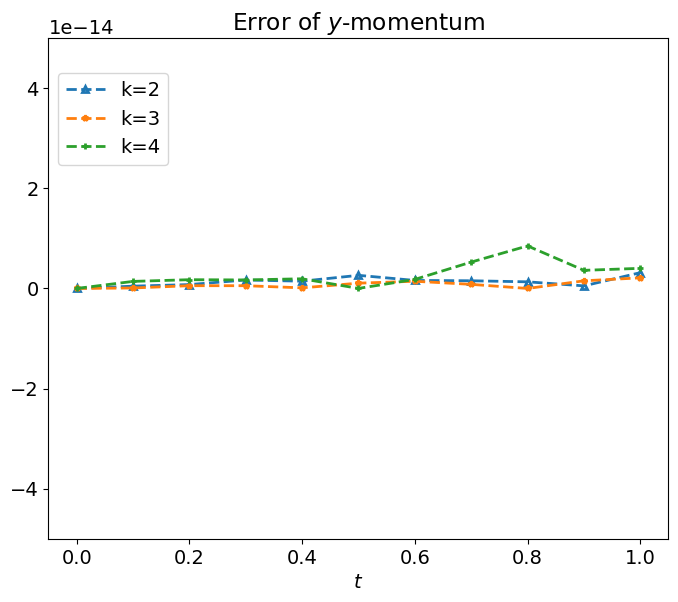}}
\caption{\small Errors of mass, energy, and momentum $($Example \ref{Example3}$)$.}\label{fig3-4}
\end{figure}

\begin{figure}[htp!] 
    \centering
    {\includegraphics[trim = .1cm .1cm .1cm .1cm, clip=true,width=0.5\textwidth,height=0.35\textwidth]{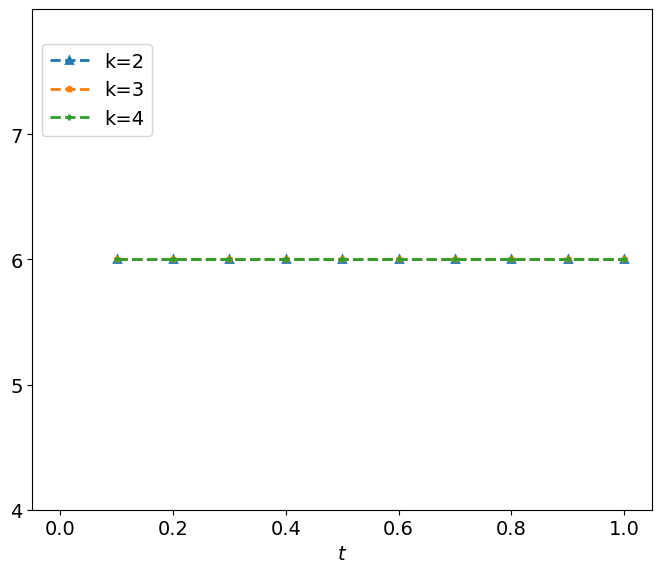}}
   \caption{\small Number of iterations at each time level $($Example \ref{Example3}$)$.}
   \label{fig33-5}
\end{figure}

Figure~\ref{fig3-1} demonstrates that the temporal discretization errors exhibit an $O(\tau^{k+1})$ convergence rate. The spatial mesh size is fixed at \( h = 5/6 \) with \( p = 1 \). In this setting, the spatial discretization error is negligible, as the reference solution is computed using the same mesh size \( h = 5/6 \) and a sufficiently small time step \( \tau = 1/100 \). This choice ensures that the temporal error in the reference solution is negligible compared to the errors corresponding to the time step sizes used in Figure~\ref{fig3-1}.
Figure \ref{fig3-2} shows that the errors from the spatial discretizations are $O(h^p)$, where the time stepsize is chosen to be 
$\tau = 1/100$ with $k = 1$. These numerical results illustrate the high-order convergence of the proposed method.

\begin{figure}[htbp]
		\centering
		\subfigure[$|u|$ at $t=0$]{
		\begin{minipage}[t]{0.47\textwidth}
				\centering
				\includegraphics[width=2.0in]{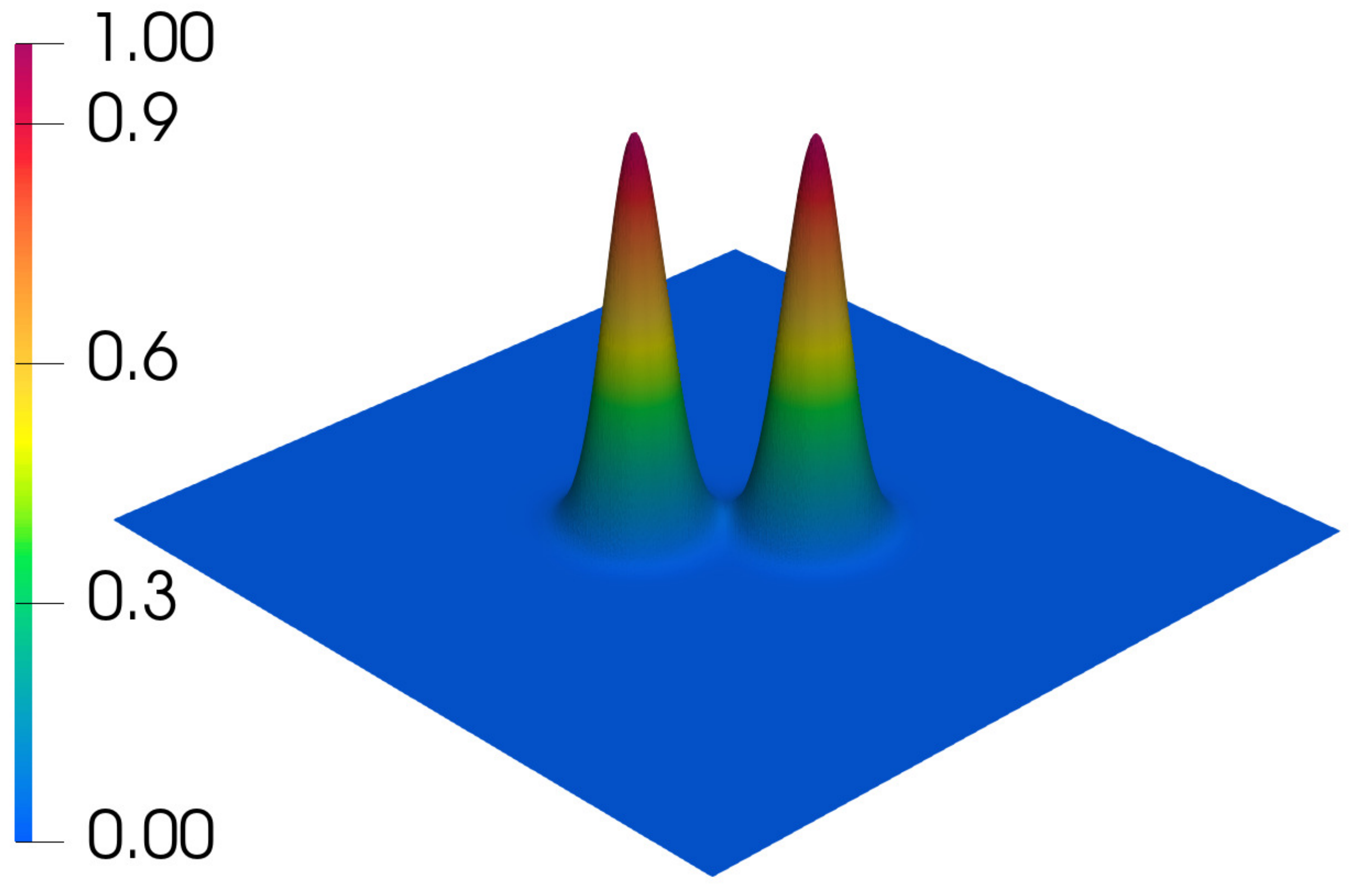} 
		\end{minipage}} 
		\subfigure[$|u|$ at $t=0.5$]{
		\begin{minipage}[t]{0.47\textwidth}
				\centering
				\includegraphics[width=2.0in]{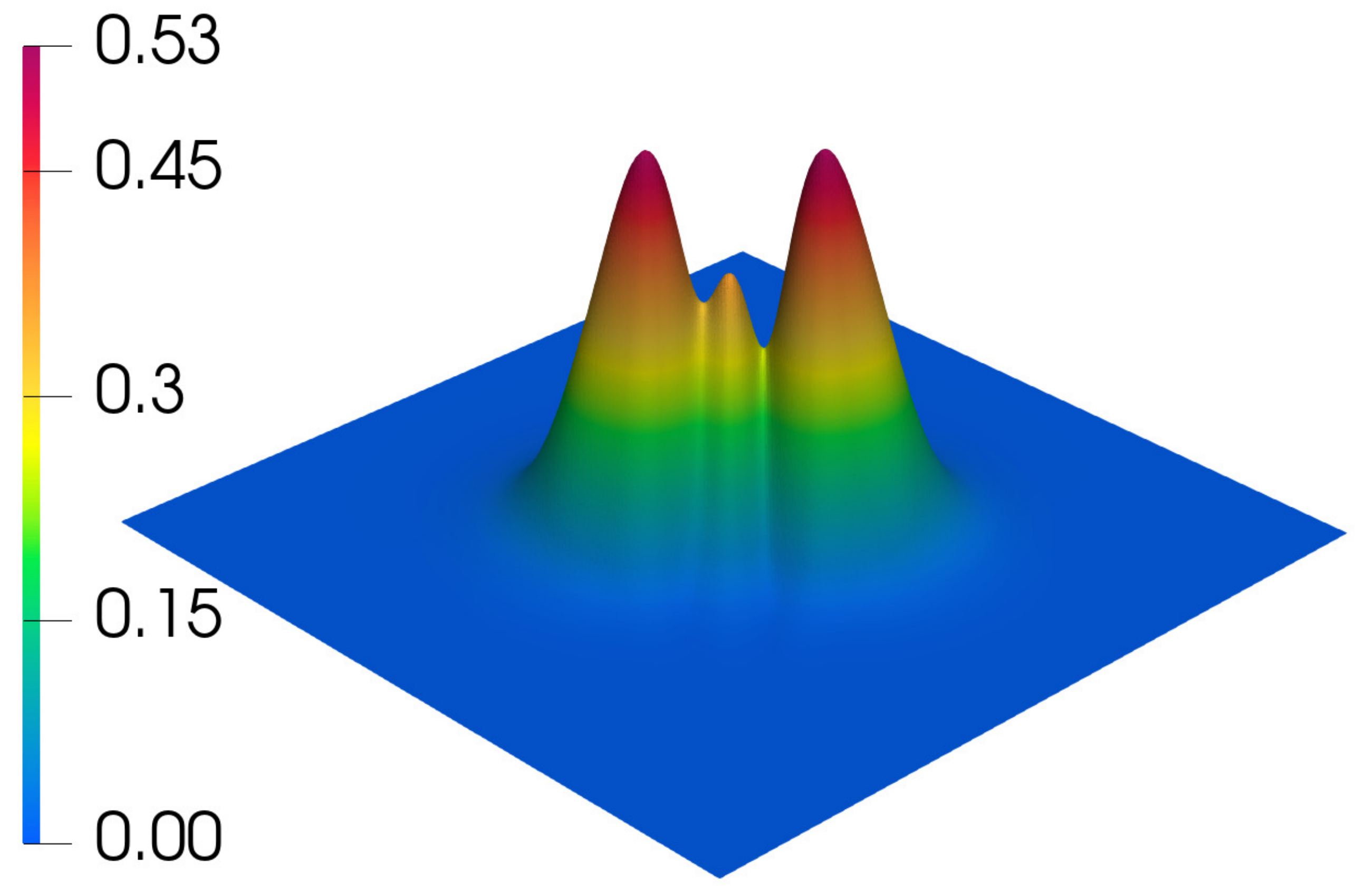} 
		\end{minipage}}
		
		\subfigure[$|u|$ at $t=1$]{
		\begin{minipage}[t]{0.47\textwidth}
				\centering
				\includegraphics[width=2.0in]{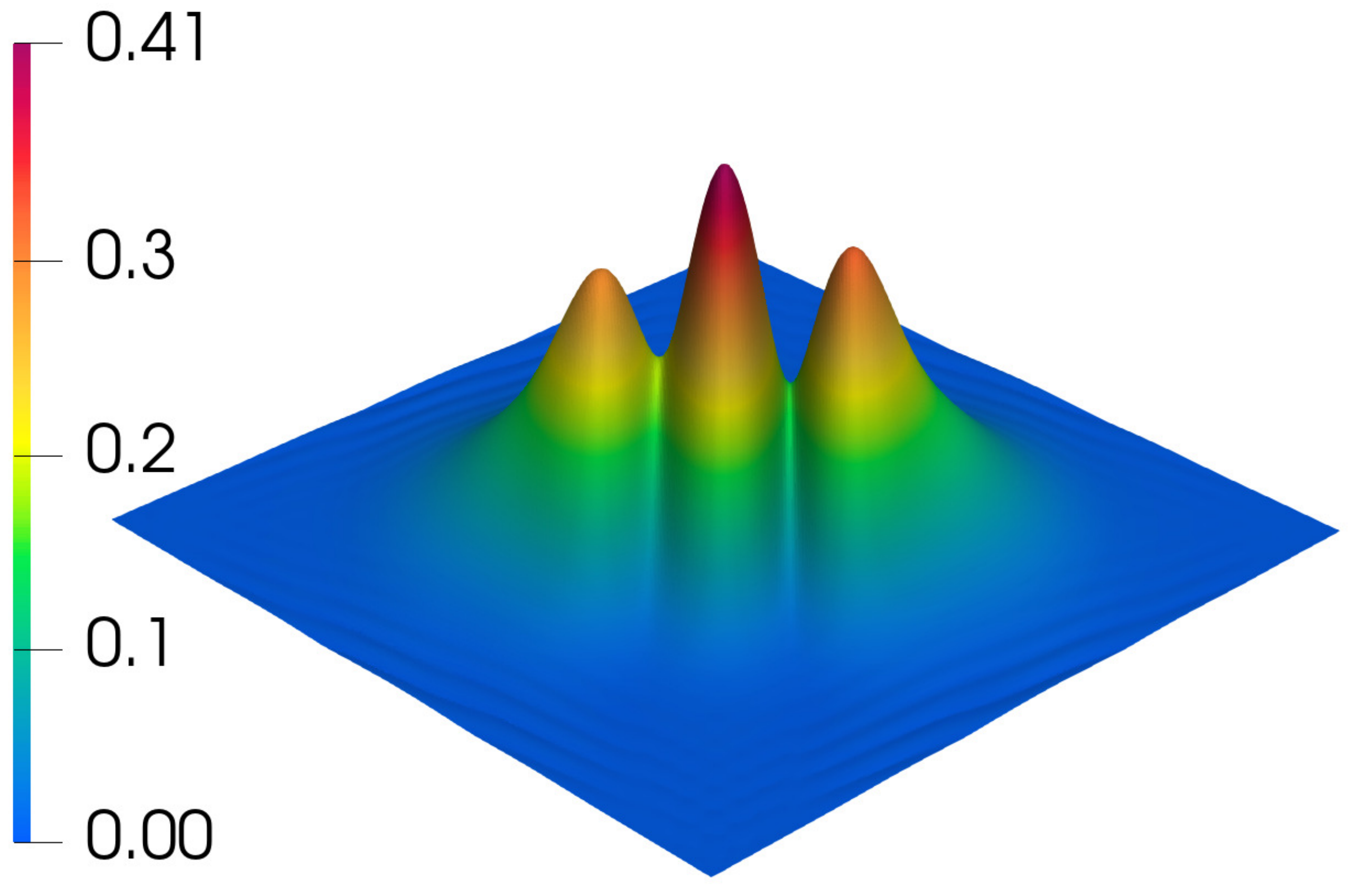} 
		\end{minipage}}
		\subfigure[$|u|$ at $t$=1.5]{
			\begin{minipage}[t]{0.47\textwidth}
				\centering
				\includegraphics[width=2.0in]{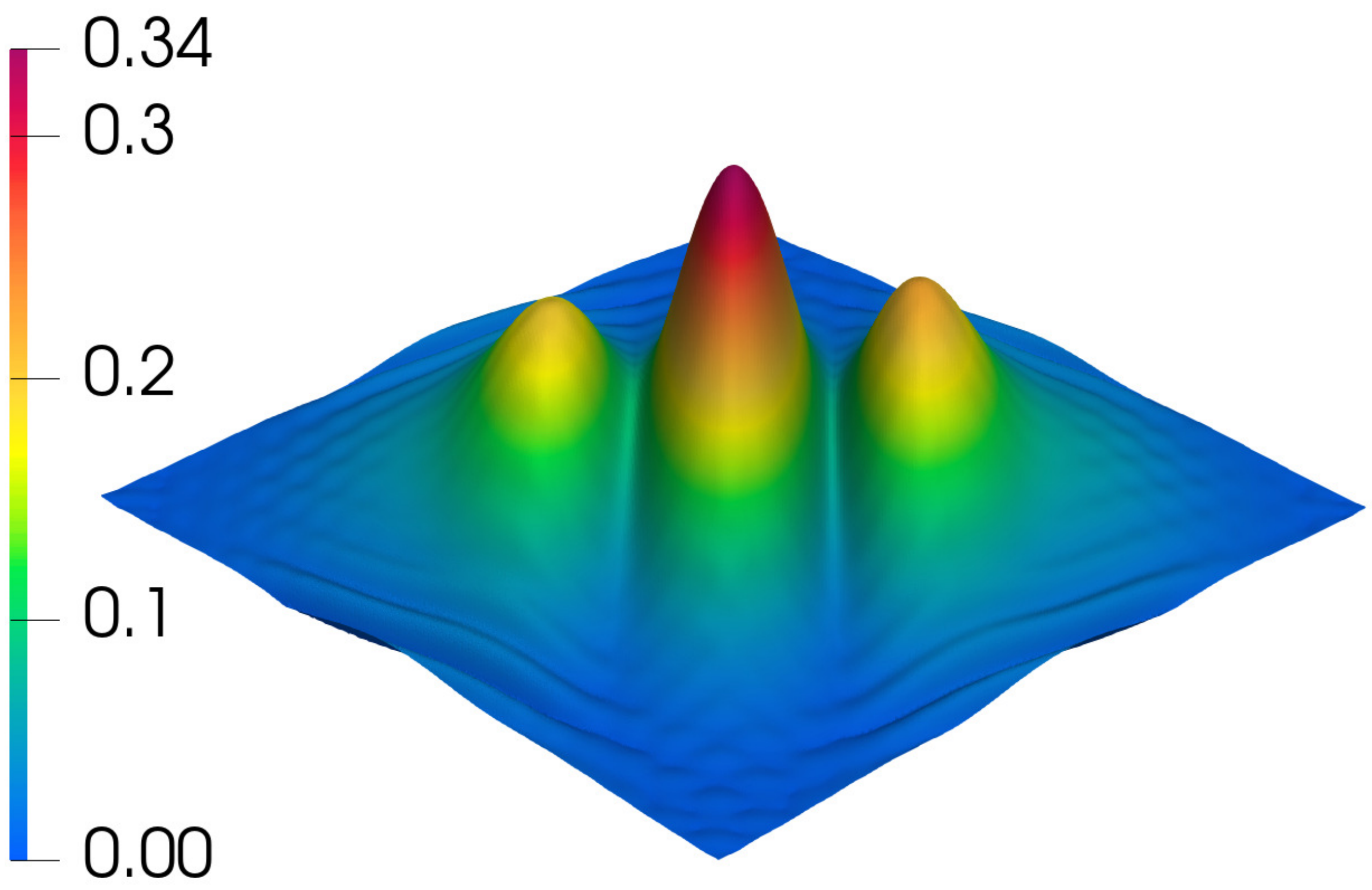}
		\end{minipage}}
\caption{Time evolutions of solitons for Example \ref{Example3}.  } 
		\label{fig3-6}
\end{figure}
\begin{figure}[htbp]
		\subfigure[$|u|$ at $t=0$]{
		\begin{minipage}[t]{0.47\textwidth}
				\centering
				\includegraphics[width=1.8in]{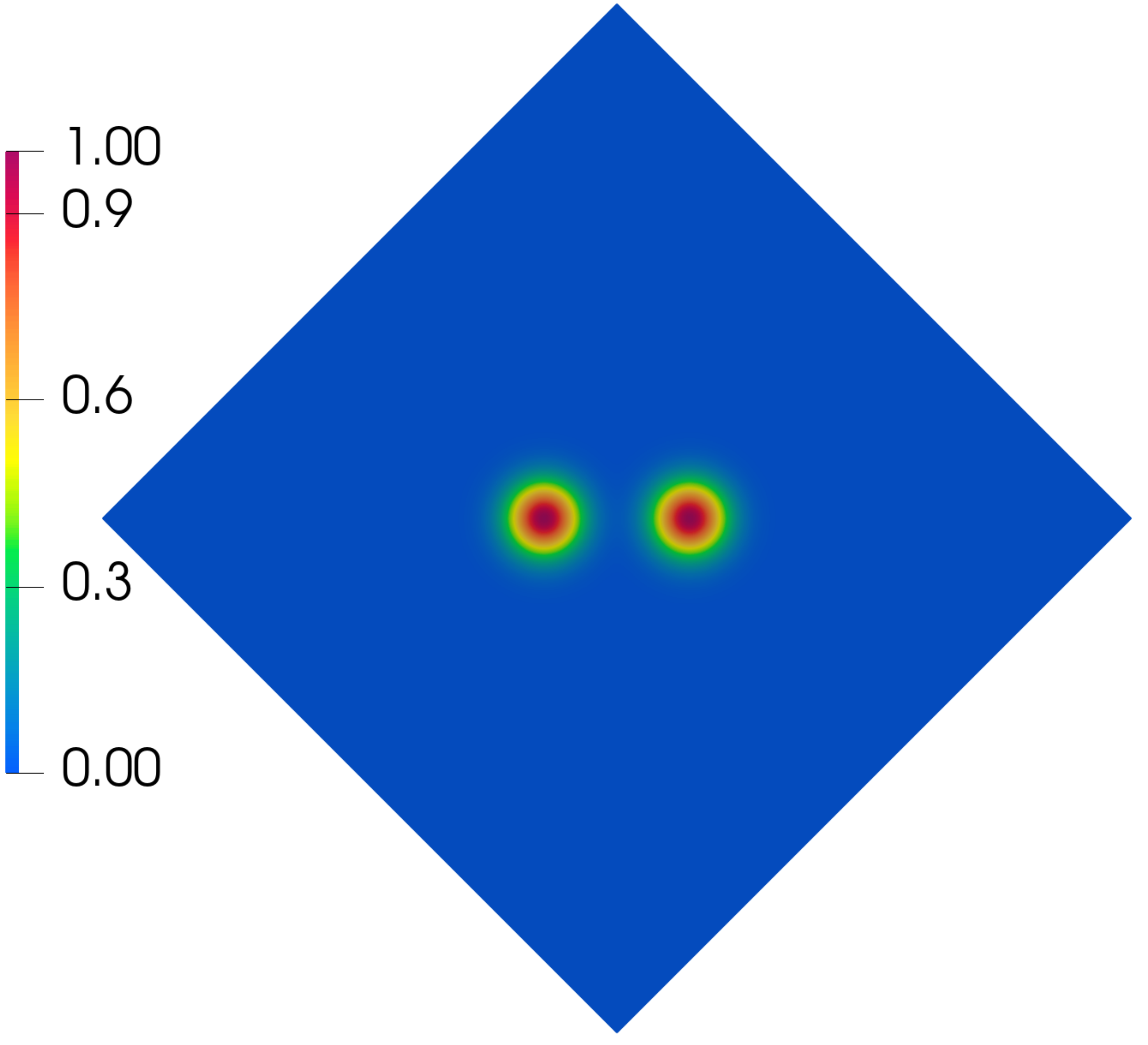} 
		\end{minipage}}
		\subfigure[$|u|$ at $t=0.5$]{
		\begin{minipage}[t]{0.47\textwidth}
				\centering
				\includegraphics[width=1.8in]{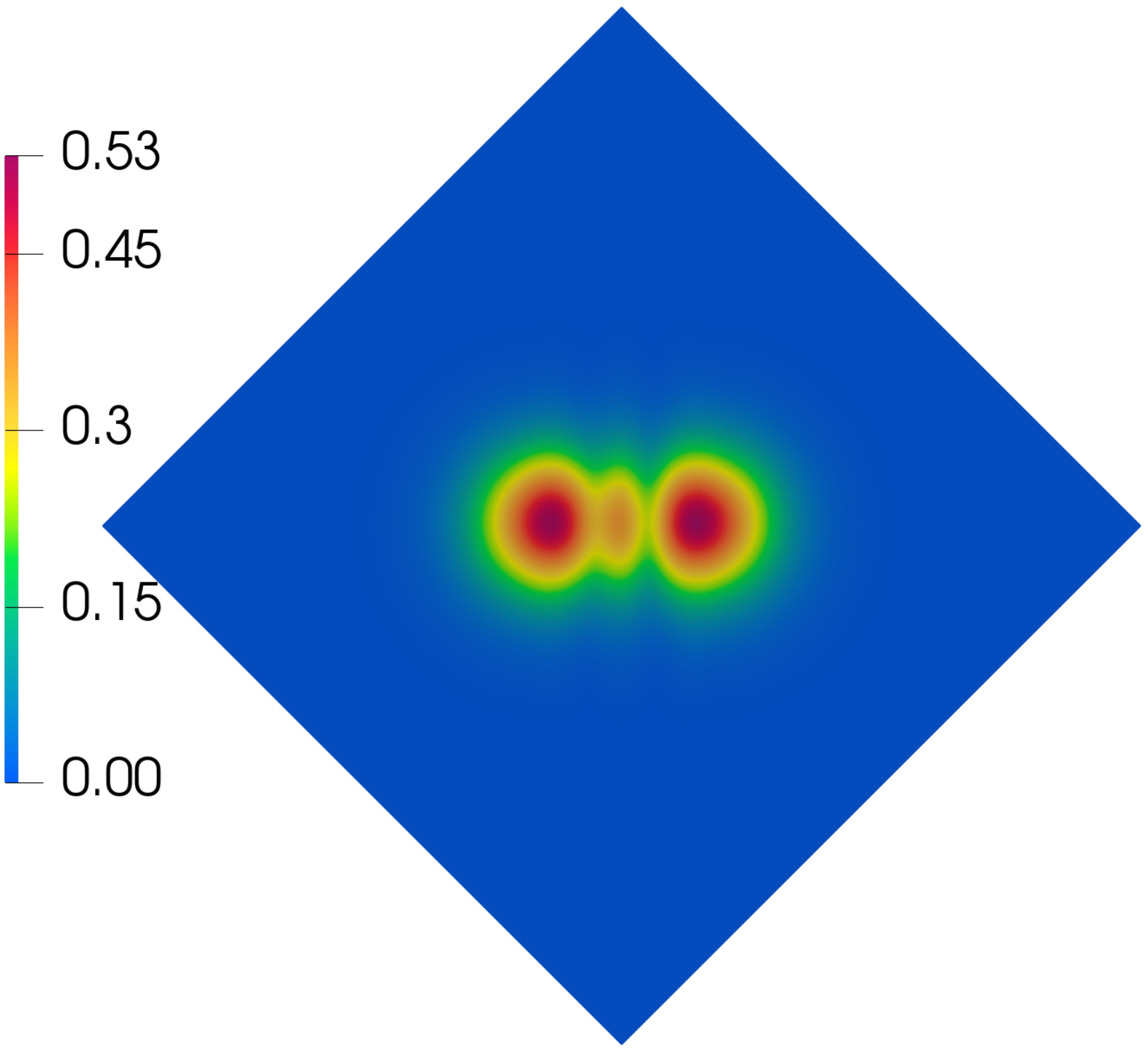} 
		\end{minipage}}
		
		\subfigure[$|u|$ at $t=1$]{
		\begin{minipage}[t]{0.47\textwidth}
				\centering
				\includegraphics[width=1.8in]{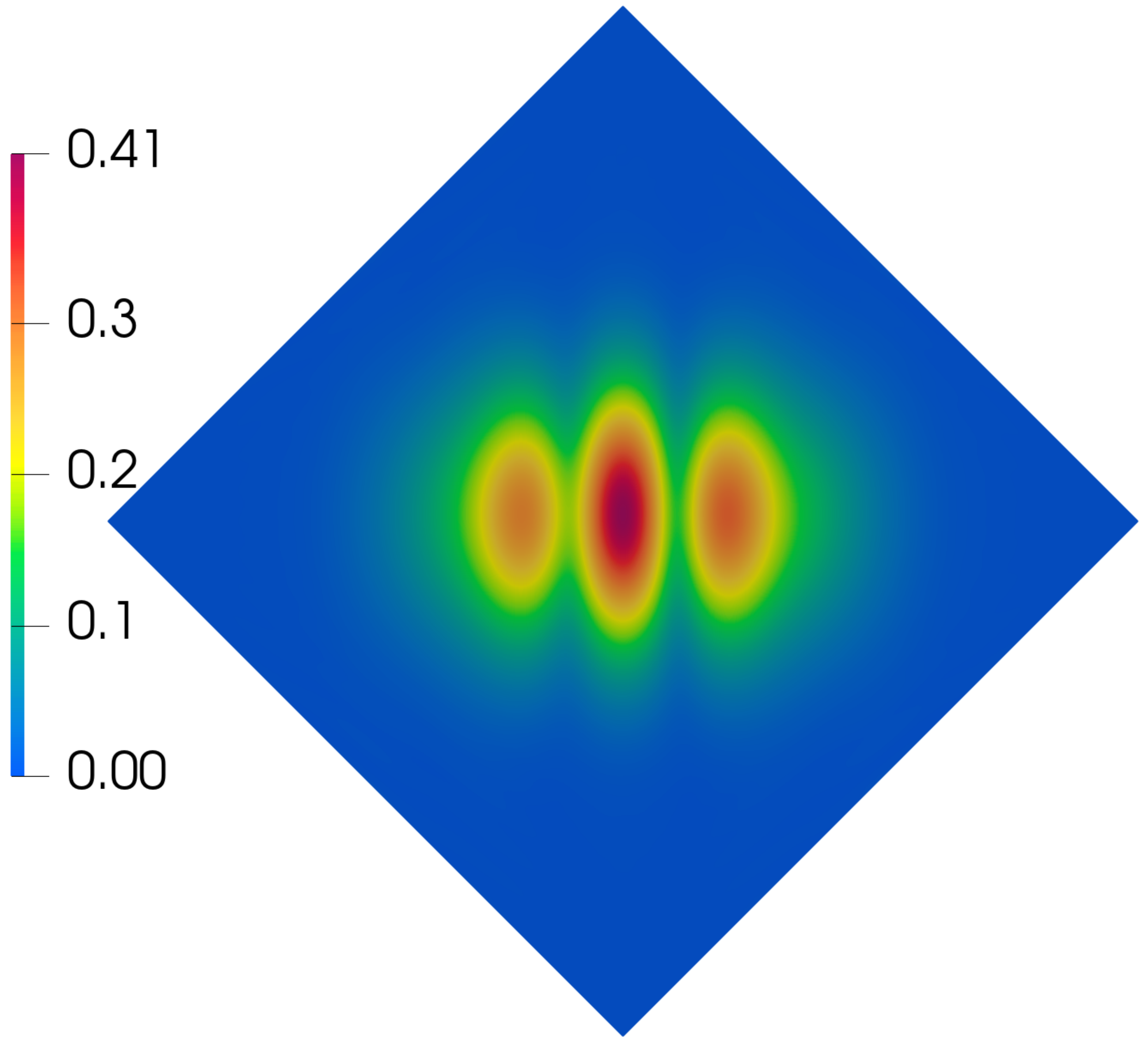} 
		\end{minipage}}
		\subfigure[$|u|$ at $t$=1.5]{
			\begin{minipage}[t]{0.47\textwidth}
				\centering
				\includegraphics[width=1.8in]{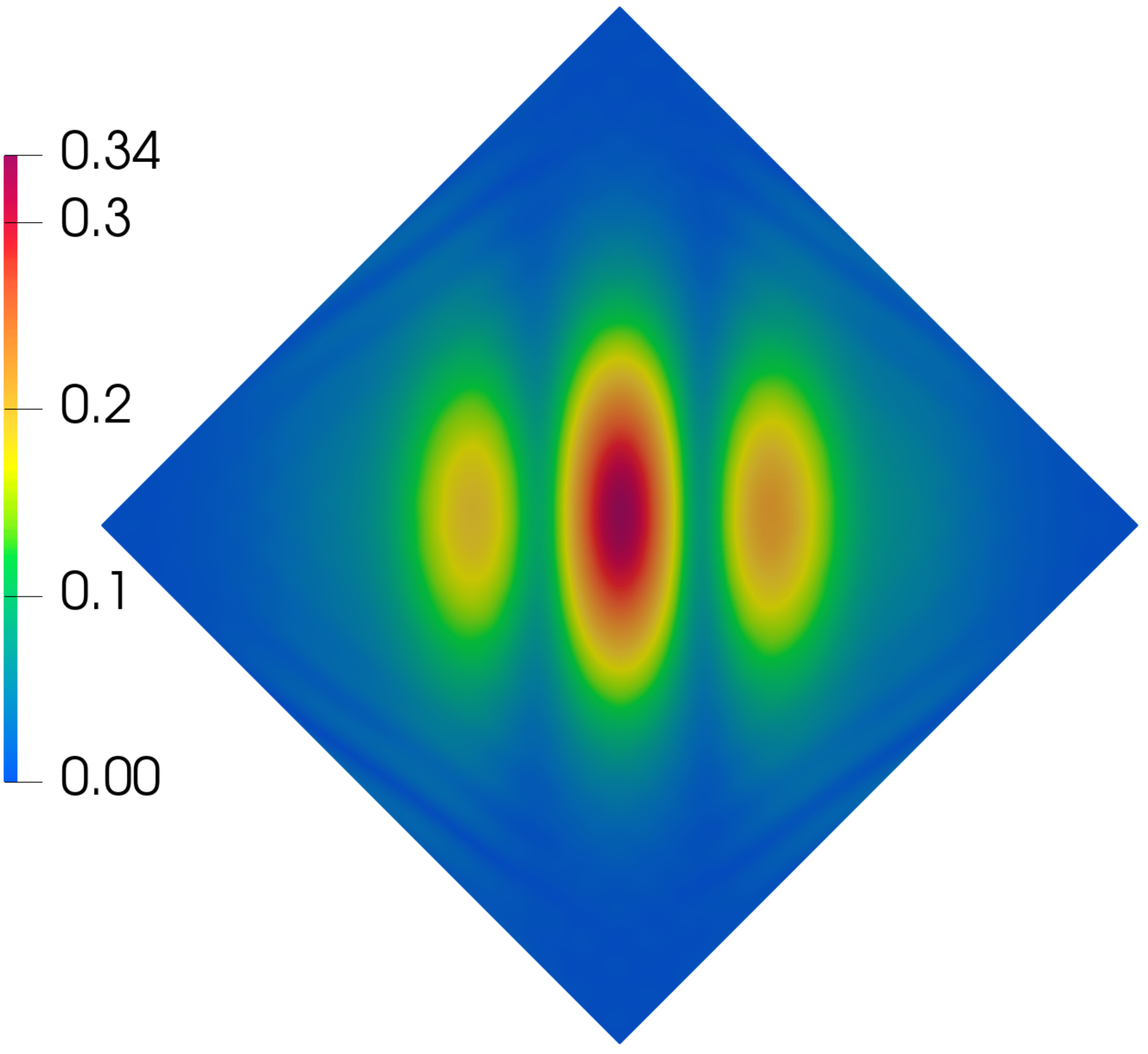}
		\end{minipage}}
\caption{Time evolutions of solitons for Example \ref{Example3}.  } 
		\label{fig3-7}
\end{figure}

The conservation of mass, energy, and momentum is presented in Figure \ref{fig3-3} with $p = 1$, $k = 2, 3, 4$, $T = 1$, 
$\tau = 1/10$, and $h = 1/2$. The errors in conserving mass, energy, and momentum are about $10^{-12}$, $10^{-12},$ 
and $10^{-14}$, respectively, as shown in Figure \ref{fig3-4}. These errors are significantly smaller than the errors 
of the numerical solutions (about $10^{-1}$ as shown in Figure \ref{fig3-2}). This shows the effectiveness of the 
proposed method in conserving the mass, energy, and momentum of the NLS equation. The numbers of iterations 
at each time level are shown in Figure \ref{fig33-5}; $6$ Newton iterations are needed at every time level. 
This is not large and it is acceptable in exchange of the conservation of mass, energy, and momentum of the numerical solutions. 

The evolution of the two-dimensional solitons is presented in Figures \ref{fig3-6}--\ref{fig3-7} for the  numerical solutions with $p = 3$, $k = 2$, 
$\tau=1/50$, and $h = 1/5$. Figure \ref{fig3-6}(a) shows that the initial state of the solution is composed of two peaks. As time increases,
 the two peaks radiate and begin to collide with each other, as shown in Figures \ref{fig3-6}(b) and \ref{fig3-7}(b). Later, the collision of solitons 
 leads to the creation of a new peak at the center of the domain, as shown in Figures \ref{fig3-6}(c) and \ref{fig3-7}(c). As time increases, 
 the amplitude of the peak at the center of the domain becomes bigger, while the amplitudes of the other two peaks become smaller; 
 see Figures \ref{fig3-6}(d) and \ref{fig3-7}(d). 

\end{example}

\begin{example}[Simulation of two-dimensional NLS equation with singular solutions]\label{Example5} \upshape
We investigate the blow-up behavior of the solution to the following NLS equation with periodic boundary condition:
  \begin{equation}\label{ex5.3-1}
      \begin{aligned}
          \i\partial_t u +\Delta u + 15|u|^2u &= 0 \quad &&\text{in } \Omega \times (0, T], \\
          u|_{t=0} &= 5\sin(2\pi x)\sin(2\pi y) \quad &&\text{in } \Omega ,
      \end{aligned}
  \end{equation}
  where \( \Omega \) is a rectangular domain with periodicity in \( x \pm y \):
  \[
      \Omega = \{ (x, y) \in \mathbb{R}^2: |x| + |y| \leqslant \sqrt{2}L \}, \quad L = 1.
  \]
  
  In this example, we set \( p = 1 \) and \( k = 1 \). The spatial mesh size is \( h = 2\sqrt{2}/100 \), and the time step \( \tau \) is dynamically adjusted based on the \( L^\infty \)-norm of \( u \) at the previous time level:
  \begin{align}\label{adaptive_time}
      \tau_n = \frac{10^{-4}}{\|u(t_{n-1})\|_{L^\infty(\Omega)}} \quad \text{for } n\geqslant 1.
  \end{align}
  Given that the initial energy is approximately $  E_h(0) = -1.477 \times 10^2$,  computed using the specified mesh size and time step size, its negativity implies that the solution is expected to blow up in finite time \cite{Debussche2002}. In particular, blow-up is detected by the rapid growth of the norms $\|u\|_{L^\infty(\Omega)}$, $ \|u\|_{W^{1,\infty}(\Omega)}$ and  $ \|\partial_t u\|_{L^\infty(\Omega)}$.
  
Our simulations show that the blow-up occurs at \(t_{\text{blow-up}} \approx 0.0007755\), at which \(\|u\|_{W^{1,\infty}(\Omega)} \approx 1200\) and \(\|\partial_t u\|_{L^\infty(\Omega)} \approx 2 \times 10^7\); see Figure \ref{fig5-1}. As illustrated in Figure \ref{fig5-2}, both the \( L^\infty \) and \( W^{1,\infty} \) norms of \( u \) experience steep growth as \( t \) approaches \( t_{\text{blow-up}} \). 


  \begin{figure}[htbp]
      \centering
      \subfigure[$|u|$ at $t=0$]{\includegraphics[width=2.3in]{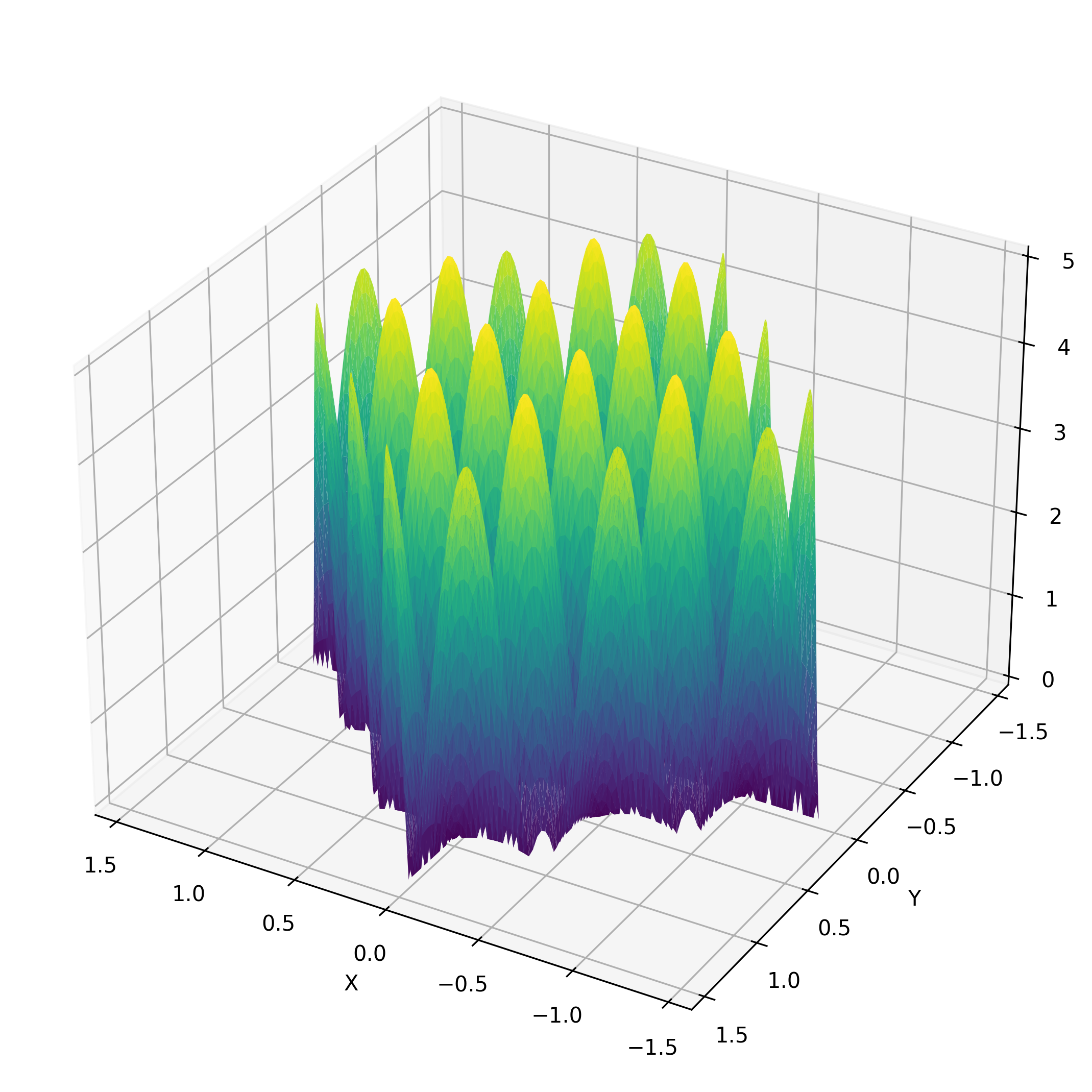}}
      \subfigure[$|u|$ at $t=t_{\text{blow-up}}$]{\includegraphics[width=2.3in]{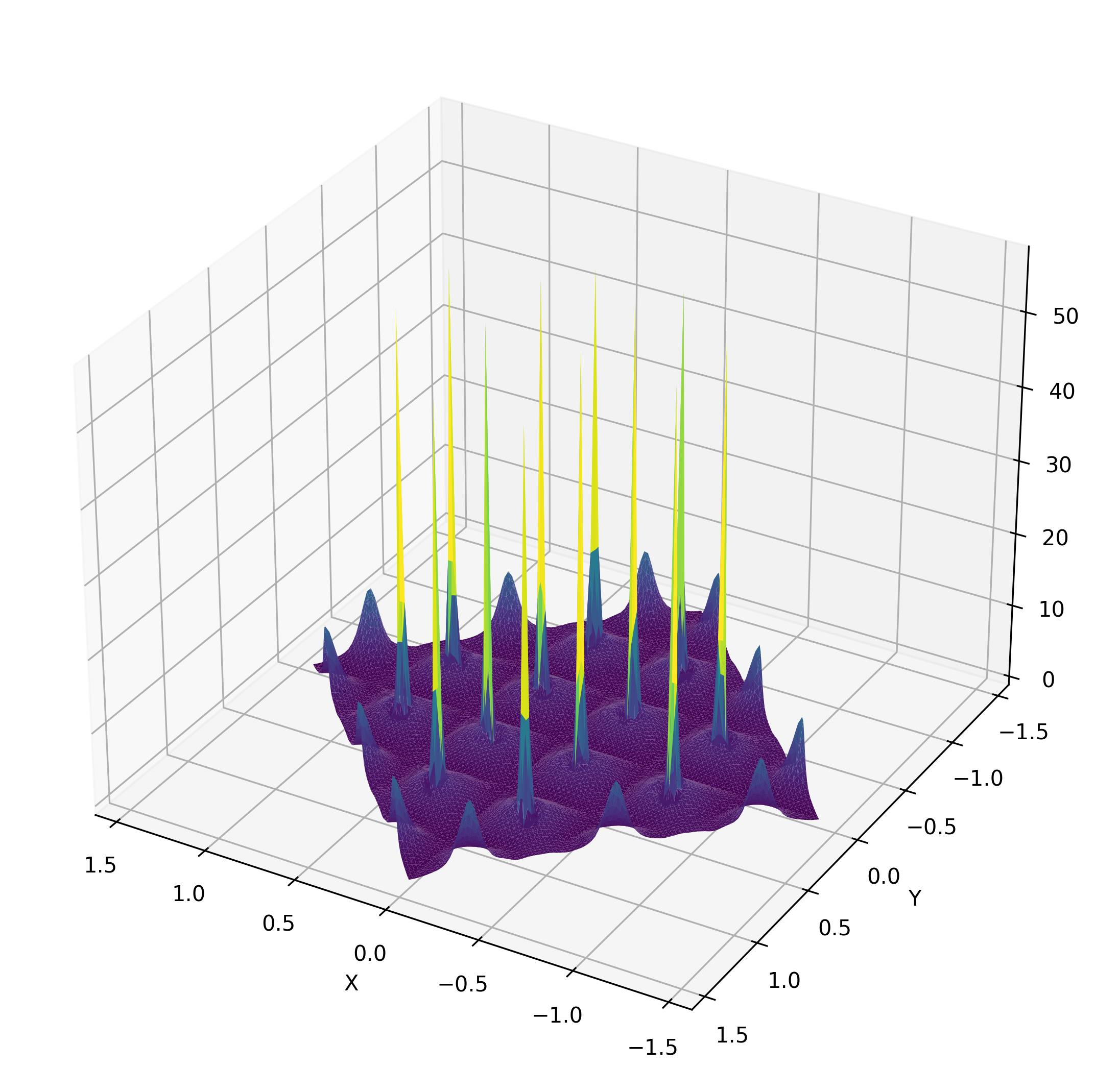}}
      \subfigure[$|\nabla u|$ at $t=0$]{\includegraphics[width=2.3in]{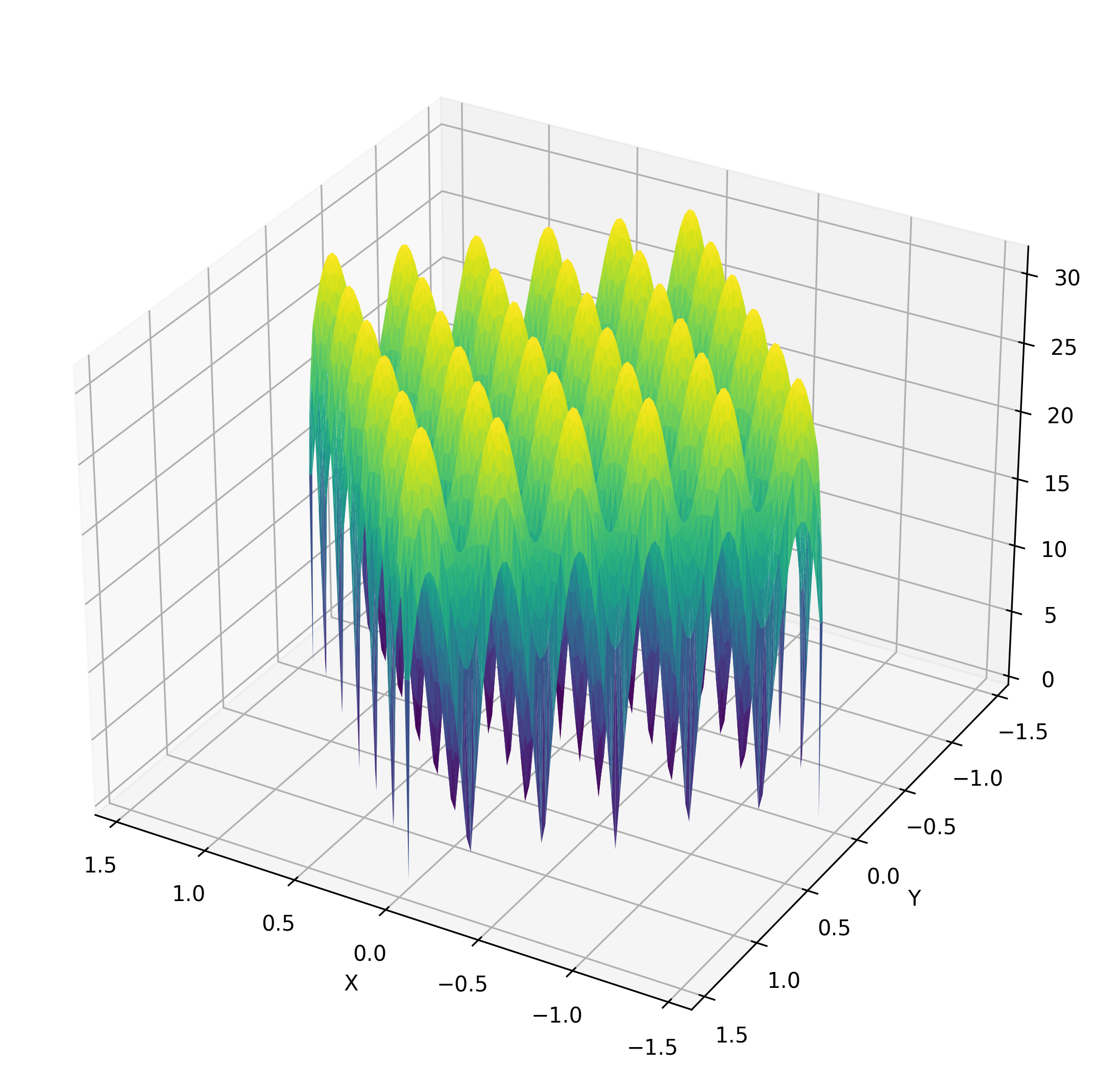}}
      \subfigure[$|\nabla u|$ at $t=t_{\text{blow-up}}$]{\includegraphics[width=2.3in]{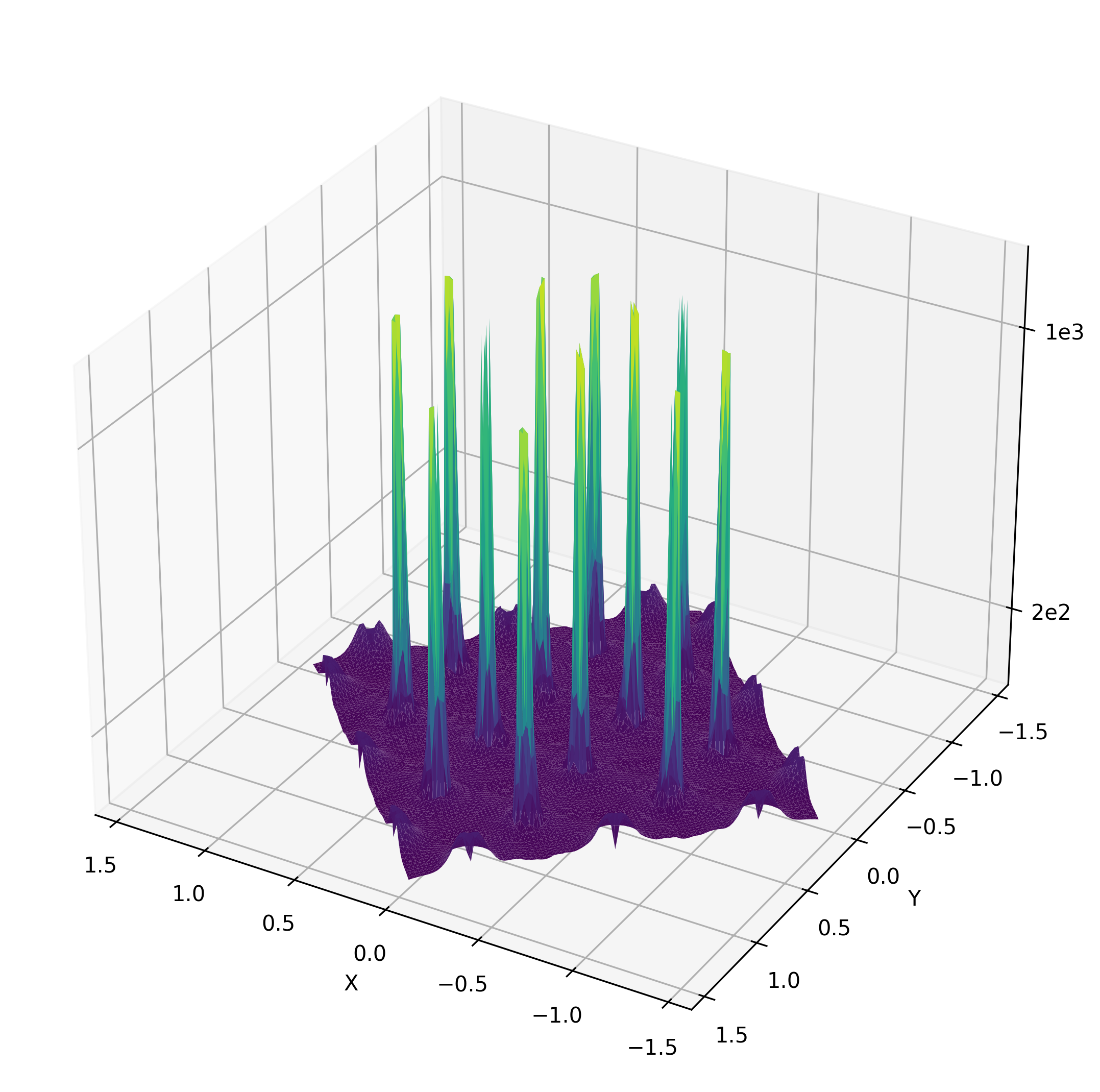}}
      \caption{Graph of $|u|$ and $|\nabla u|$ at $t=0$ and $t_{\text{blow-up}}$ (Example \ref{Example5}).}
      \label{fig5-1}
  \end{figure}

  \begin{figure}[htbp]
      \centering
      \subfigure{\includegraphics[width=2.8in]{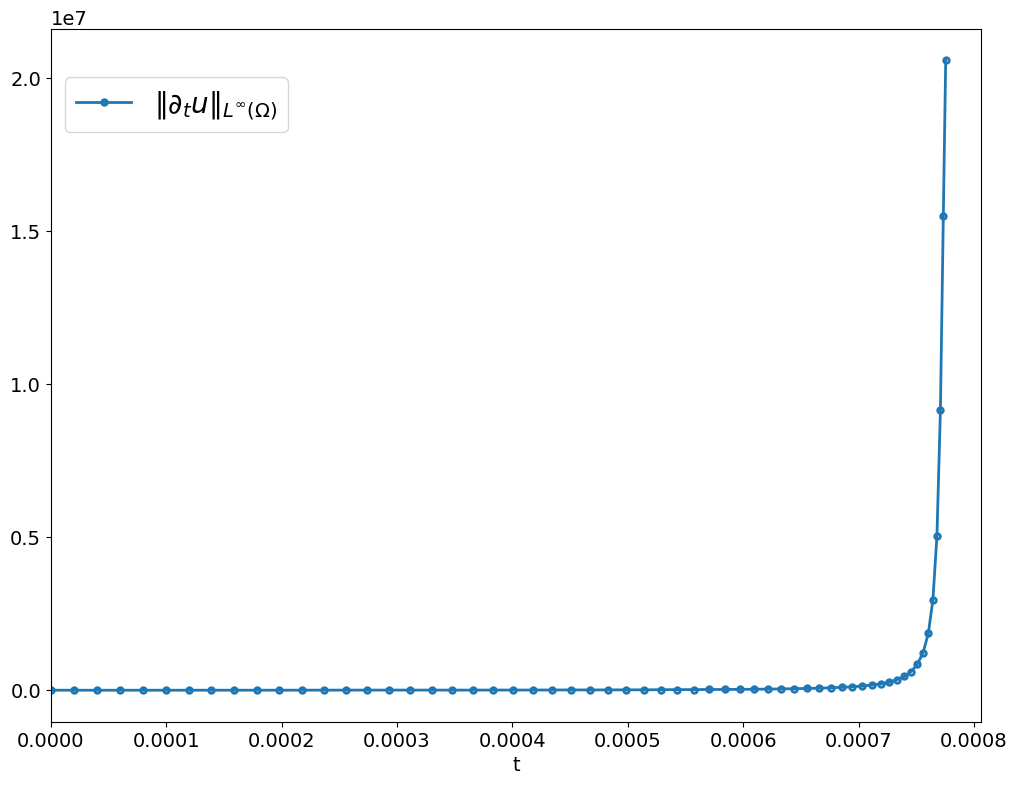}}
      \subfigure{\includegraphics[width=2.8in]{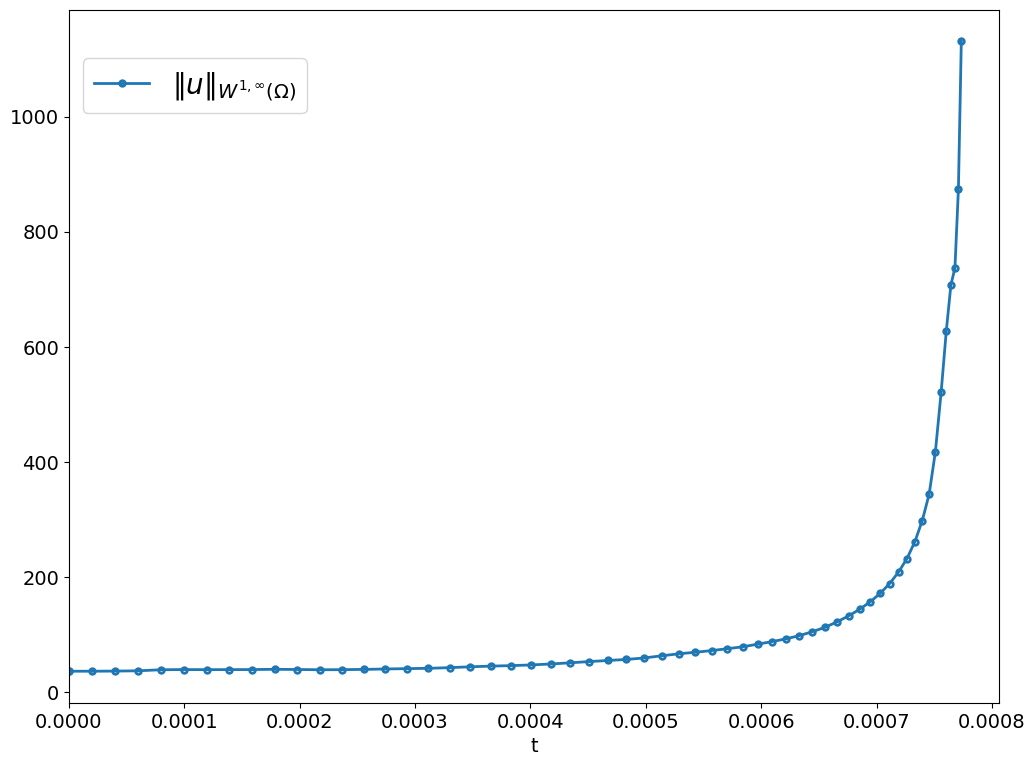}}
      \caption{Evolutions of $\|\partial_t u\|_{L^\infty(\Omega)} $ and $\|u\|_{W^{1,\infty}(\Omega)} $, up to blow-up time (Example \ref{Example5}).}
      \label{fig5-2}
  \end{figure}
  
The proposed nonlinearly implicit scheme conserves energy exactly; however, the incomplete Newton iteration introduces errors that become more pronounced as the solution approaches blow-up. To ensure that energy loss remains below a certain threshold, we have implemented the following stopping criteria throughout:
  \begin{subequations}\label{stopping-criteria}
    \begin{align}
    \| u_h^{(\ell)}(t_n) - u_h^{(\ell-1)}(t_n) \|_{H^1(\Omega)} &< 10^{-9},\label{stopping-a}\\
    |E(u_h^n) - E(u_h^{n-1})|& < 10^{-8}, \label{stopping-b}
    \end{align}
  \end{subequations}
With these stopping criteria, the mass, energy and momentum remain conserved within errors of \(10^{-11}\), \(10^{-8}\), and \(10^{-12}\), respectively, up to the blow up time, as illustrated in Figure~\ref{fig5-3}. In particular, the energy conservation is maintained to within \(10^{-8}\), which aligns with the termination condition in Newton's iteration. Therefore, the second criterion is crucial to ensure that the energy error remains sufficiently small (i.e., on the order of \(10^{-8}\)) relative to the magnitude of the solution. Figure \ref{fig5-4} illustrates the number of iterations performed per time step. In particular, as \(t\) approaches the blow-up time, the number of iterations also increases to guarantee that the energy loss is within the tolerance defined in \eqref{stopping-b}.

  \begin{figure}[htp!] 
      \centering
      \includegraphics[width=0.44\textwidth,height=0.36\textwidth]{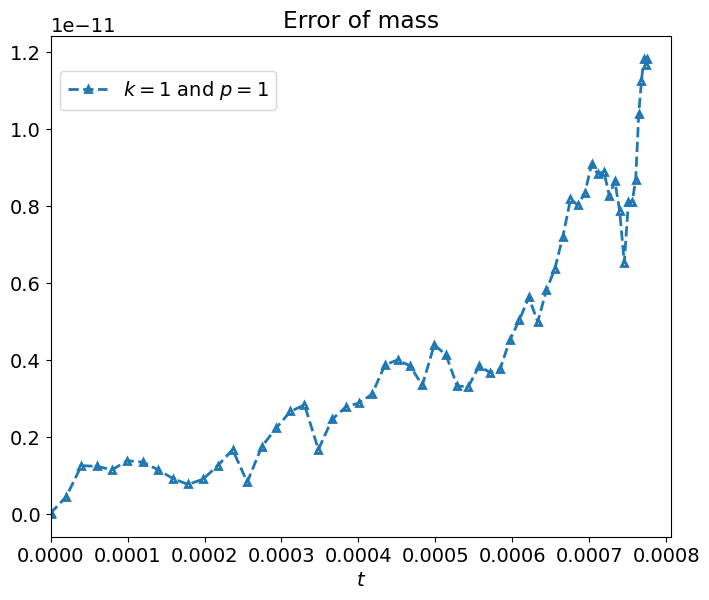}
      \includegraphics[width=0.44\textwidth,height=0.36\textwidth]{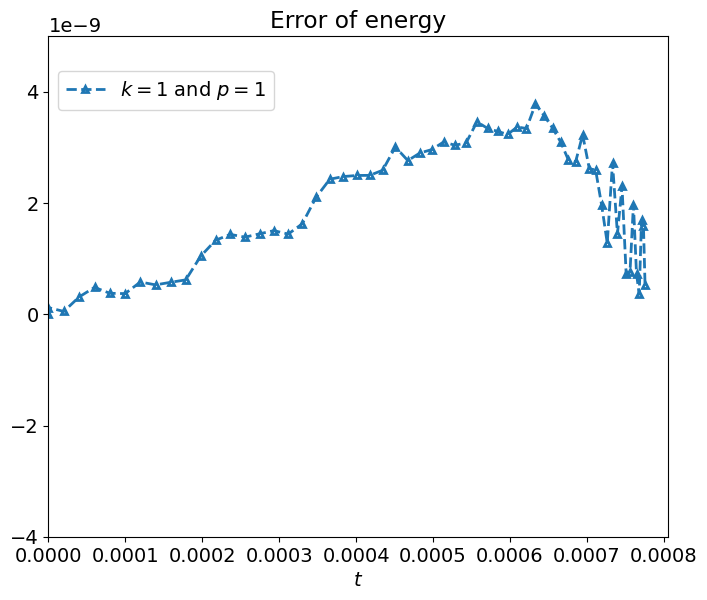}
      \includegraphics[width=0.44\textwidth,height=0.36\textwidth]{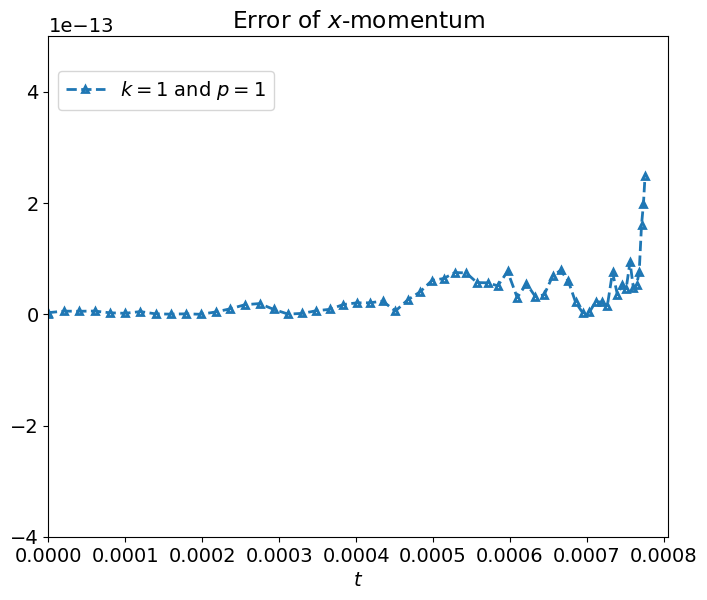}
      \includegraphics[width=0.44\textwidth,height=0.36\textwidth]{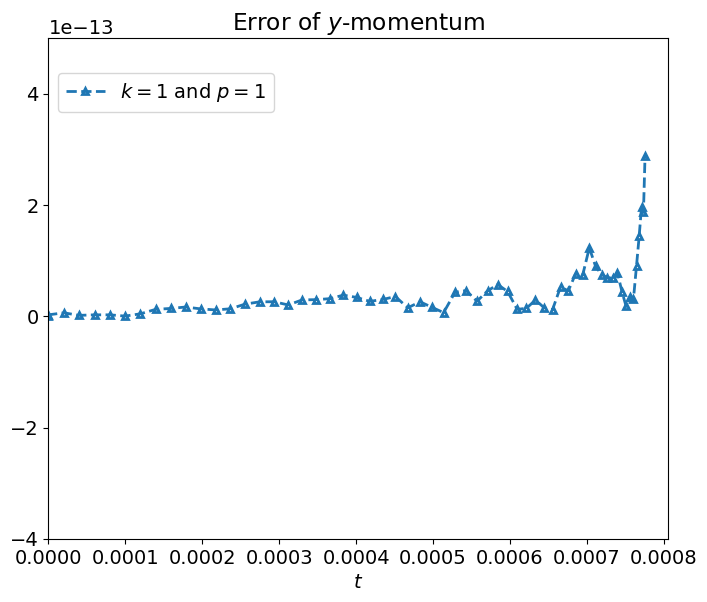}
      \caption{Errors of mass, energy, and momentum (Example \ref{Example5}).}
      \label{fig5-3}
  \end{figure}

  \begin{figure}[htp!] 
      \centering
      \includegraphics[width=0.45\textwidth,height=0.40\textwidth]{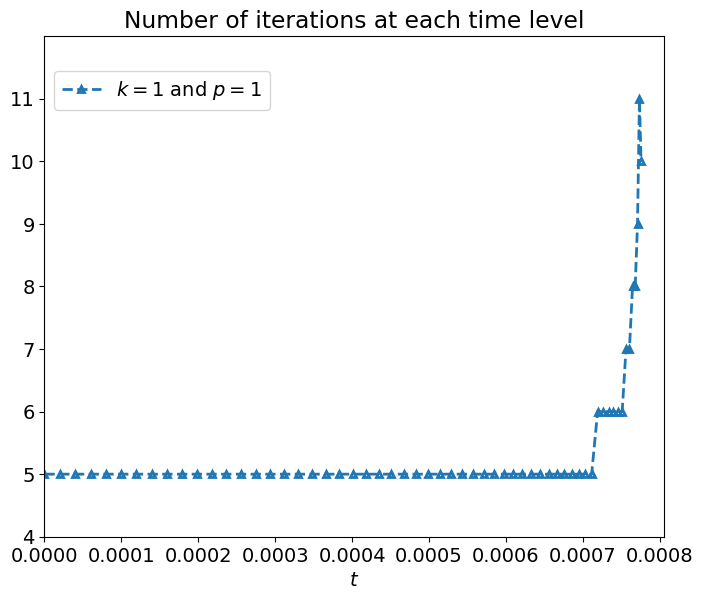}
      \caption{Number of iterations per time level (Example \ref{Example5}).}
      \label{fig5-4}
  \end{figure}
  \end{example}

\section{Conclusion}%
\label{conclusion}
We have introduced a novel formulation along with an associated space-time finite element method (FEM) for the NLS equation. 
We have demonstrated that the proposed algorithm conserves mass, energy, and momentum of the NLS equation at the discrete level, 
and we have designed a semi-Newton iteration for the nonlinear system associated with the numerical method. Through extensive 
numerical examples, we have shown that the proposed method achieves high-order convergence in approximating solutions 
of the NLS equation and effectively conserves mass, energy, and momentum for the simulation of various solitons of the NLS equations.

\section*{Acknowledgement}
The work of B. Li is supported in part by the National Natural Science Foundation of China (Project No. 12231003), 
the Hong Kong Research Grants Council (GRF Project No. PolyU15306123), and an internal grant of The Hong Kong Polytechnic University (Project ID: P0045404), and a Collaborative Research Project at The Hong Kong Polytechnic University (Project ID: P0038443). The work of R. Tang is supported by the AMSS-PolyU Joint Laboratory and Hong Kong Research Grants Council (GRF Project No. PolyU15303022). 
The work of H. Zhang is supported by the National Natural Science Foundation of China (Project Nos. 12120101001,12371447,12171284), the Natural Science Foundation of Shandong Province (Project Nos. ZR2021ZD03), and the Hong Kong Research Grants Council (GRF Project No. PolyU15301321).



\end{document}